\documentclass[12pt,psamsfonts]{article}
\newcommand{\vrpf}[2]{\begin{picture}(#1,1)
\put(0,1){\vector(1,0){#1}}
\put(0,2.5){\makebox(#1,0)[b]{$#2$}}
\end{picture}}

\newcommand{\rupfl}[1]{\begin{picture}(12,12)
\put(0,12){\vector(1,-1){12}}
\put(8,0){\makebox(0,12)[l]{$#1$}}
\end{picture}}

\newcommand{\ropfr}[1]{\begin{picture}(12,12)
\put(0,0){\vector(1,1){12}}
\put(4,0){\makebox(0,12)[r]{$#1$}}
\end{picture}}

\usepackage{amsmath}
\usepackage{amssymb}
\newtheorem{theo}{Theorem}[section]
\newtheorem{remarkk}[theo]{Remark}
\newenvironment{rem}{\begin{remarkk}\rm}{\end{remarkk}}
\newtheorem{definition}[theo]{Definition}
\newenvironment{defi}{\begin{definition}\rm}{\end{definition}}
\newtheorem{prop}[theo] {Proposition}
\newtheorem{cor}[theo]{Corollary}
\newtheorem{lemma}[theo]{Lemma}
\newtheorem{example}[theo]{Example}

\begin{document}
\title{Canonical rings of surfaces whose canonical system has base points}
\author{I. C. Bauer, F. Catanese, R. Pignatelli\thanks{The present cooperation took place in the realm of the
  DFG-Forschungsschwerpunkt ``Globale methoden in der komplexen
  Geometrie'' and of the EAGER project.}
}

\date{\today}
\maketitle

\tableofcontents

\section{Introduction}
In Enriques' book on algebraic surfaces (\cite{Enr}), culminating
a research's lifespan of over 50 years, much emphasis was set
on the effective construction of surfaces, for instance of  surfaces with
$p_g =4$ and whose canonical map is a birational map onto a singular
surface
$\Sigma$ in  $\mathbb{P}^3$.

The problem of the effective construction of such surfaces for
the first open case $K^2 = 7$ has attracted the attention of several
mathematicians, and special constructions have been obtained by
Enriques (\cite{Enr}),
Franchetta (\cite{Fran}), Maxwell (\cite{Max}), Kodaira (\cite{Kod}).
Until  Ciliberto
(\cite{Cil1}) was able to construct an irreducible  Zariski open set
 of the  moduli space of
(minimal algebraic) surfaces with
$K^2 = 7$,
$p_g =4$, and constituted  by surfaces  with a birational canonical
morphism
whose image $\Sigma$ has ordinary singularities.

Later on,  through work of the first two named authors and of Zucconi
(\cite{Ba}, \cite{CatBu}, \cite{Zuc}), the complete classification of
surfaces with $K^2 =
7$,
$p_g =4$ was achieved, and it was shown in \cite{Ba} that the moduli
space consists of three ireducible components (two of them
consist of surfaces with non birational canonical map). But, as in
the previous work of Horikawa (\cite{HorI-V},
\cite{HorQ}) who classified surfaces with  $K^2 = 5,6$, $p_g =4$, a
complete picture of
the moduli space is missing (for instance, it is still an open question
whether the moduli space for  $K^2 =7$,
$p_g =4$ has one or two connected components).

Usually, classifying surfaces with given invariants $K^2 , p_g$, is
achieved by writing a finite number of families such that every
such surface
occurs in precisely one of those families.
Each family yields a locally closed stratum of the moduli space,
and the basic question is how are these strata patched together.

Abstract deformation theory is very useful since Kuranishi's theorem
(\cite{Kur}) gives a
lower bound for the local dimension of the moduli space, thus it helps to
decide which strata are dominating a component of the moduli space.

In principle, the local structure of the moduli space
(\cite{Palamodov})
is completely described by a sequence of Massey products on the tangent
cohomology of the surface, and Horikawa clarified the structure of the moduli
space in the "easy" case of numerical 5-ics ($K^2 = 5,p_g =4$) by using the Lie
bracket
  $ H^1(S,\Theta_S) \times H^1(S,\Theta_S) \rightarrow H^2(S,\Theta_S)$.

  However, the analytic approach does not make us see concretely how do
surfaces belonging to one family deform to surfaces in another family,
and therefore Miles Reid, in the Montreal Conference of 1980 proposed to look
at the deformations of the canonical rings for numerical 5-ics
(cf.\cite{Rei0}).

His program was carried out by E. Griffin (\cite{Gri}) in this case, later
on D. Dicks found an interesting approach to the question and applied it to
the case of surfaces with $K^2 = 4, p_g =3$
(\cite{Dicks},\cite{Dicks2}).
His method was clearly exposed in the article (\cite{Rei2}) by Miles
Reid, where he set as a
challenge the problem to apply these methods to the hitherto  still partially
unexplored case of surfaces with $K^2 = 7$, $p_g=4$.

In \cite{CatBu} was given a method (of the so called quasi generic canonical
projections) allowing in principle to describe the  canonical rings of surfaces
of general type. The method works under the assumption that
the surface admits a morphism to a 3-dimensional projective
space which is a projection of the canonical map, and is birational to its
image $\Sigma$.

What happens when the canonical system has base points, in particular in our
case of surfaces with $K^2 = 7$, $p_g=4$?

Thus the first  aim of this paper is
to introduce a general method to calculate the canonical ring of minimal
surfaces
of general type whose canonical system has base points but yields a birational
canonical map. \\

We will then apply this
method  in the case of the surfaces $S$ with
$K^2 = 7$,
$p_g =4$. We will  compute the canonical ring of those minimal smooth
algebraic surfaces $S$ with
$K^2 = 7$,
$p_g =4$, whose canonical system has just one simple base point and gives a
birational map from
$S$ onto a sextic in $\mathbb{P}^3$: this is the only case, for these
values of $K^2$, $p_g$, where the canonical system has base points,
but yields a birational map. \\

What does it mean to compute a ring? As a matter of fact, using the computer
algebra program Macaulay II, we will give three different descriptions of the
above canonical rings. These presentations will allow us to deform explicitly
the canonical ring of such a minimal surface (with $K^2 = 7$, $p_g =4$ and
with birational canonical map onto a
sextic in $\mathbb{P}^3$) to the canonical ring of a surface with the same
invariants but with base point free canonical system. \\

That these deformations should exist
was already seen in
\cite{Ba}, since it was proven there that the surfaces with
$K^2 = 7$,
$p_g =4$ such that the canonical
system gives a birational map from $S$ onto a sextic in
$\mathbb{P}^3$ form an irreducible family
of dimension $35$ in the moduli space $\mathfrak{M}_{K^2 = 7, p_g =
4}$ and therefore they cannot dominate an
irreducible component of the moduli space (by Kuranishi's theorem the
dimension of
$\mathfrak{M}_{K^2 = 7, p_g =
4}$ in any point has to be at least $10\chi - 2K^2 = 36$).

Therefore it was clear from the
classification given in \cite{Ba} that this family has to be contained in
the irreducible component of the moduli space whose general point
corresponds to
a surface with base point free canonical system (obviously then with birational
   canonical morphism).

Enriques (\cite{Enr}) proposed to obtain this deformation starting by
   the surface of degree seven, union of the sextic surface (the canonical
image) together with the plane containing the double curve: in our case,
however, we see that the canonical images of  degree seven do indeed degenerate
 to the union of the sextic canonical surface  together with another plane,
namely the tacnodal plane (cf. section 3).

\bigskip
Now, although our method applies in a much more general setting, the
complexity of the computations which are
needed in every specific case grows incredibly fast. \\

We consider therefore a real challenge for our present days computer algebra
programs to make it possible to treat  surfaces with higher values of the
invariant $K^2$.\\

We would however like to remark, that all our explicit
computations are more "computer assisted
computations" than computer algebra programs. That is: it would be almost
impossible to do them without a computer
algebra program, but on the other hand there are always several steps which
have to be done by hand, because looking
carefully with a mathematical eye we can see tricks that the computer
alone cannot detect.
\\

Our paper is
organized as follows.
\\

In the first chapter we introduce under quite general conditions a naturally
defined graded subring
$\tilde{\mathcal{R}}$ of the canonical ring $\mathcal{R}$, such that
there is an exact sequence

$$
0 \longrightarrow \tilde{\mathcal{R}}_m \longrightarrow \mathcal{R}_m
\longrightarrow
H^0(\tilde{S}, \mathcal{O}_{m {\mathcal E}}) \longrightarrow
\mathbb{C}^r \longrightarrow 0,
$$

where ${\mathcal E}$ is an exceptional divisor on the surface $\tilde{S}$
obtained from $S$ blowing
up  the base points of the canonical
system. Then we introduce the ``dual'' module $M=Ext^1(\tilde{\mathcal
R}, \Gamma_*(\omega_{\mathbb P^3}))$.

The rough idea is now to calculate the subring
$\tilde{\mathcal{R}}$ (and the dual module $M$) 
using the geometry of the canonical image of
$S$.
We proceed in each degree ``enlarging''
$\tilde{\mathcal{R}}$ to  $\mathcal{R}$: we will see how the module 
$M$ provides automatically a certain number of the ``missing''
generators and relations; the few  
remaining generators and relations have to be
computed ``by hand'' by the above exact sequence.\\

In chapter $2$ we will run this program in the special case of
surfaces with $K^2 = 7$, $p_g = 4$,
whose canonical map has exactly one base point and is
birational. The main result of this section leads to the following

\begin{theo}
The canonical ring of a surface with $p_g=4$, $K^2=7$, such that
$|K|$ has one simple
base point and $\varphi_{K}$ is birational, is of the form ${\mathcal
R}:={\mathbb
C}[y_0,y_1,y_2,y_3,w_0,w_1,u]/I$ where the respective degrees  of the
generators of ${\mathcal
R}$ are $(1,1,1,1,2,2,3)$.

There exist a quadratic polynomial $Q$ and
a polynomial $B$ of degree $4$ such that $I$ is generated by the
$2
\times 2$ minors of the
matrix

$$A:=\left(
\begin{array}{cccc}
y_1 & y_2 & y_3 & w_1\\
w_0 & w_1 & Q(y_i) & u
\end{array}
\right),
$$
and by three more polynomials (of respective degrees $4,5,6$), the first one
of the form $-w_1^2+B(y_i)+\sum \mu_{ijk} y_i y_j w_k$, and the
other $2$  obtained from the first via the method of rolling (cf. section
$2$, in particular they have the form
$-w_1u+\ldots$, $-u^2+\ldots$).
\end{theo}

In chapter $3$ we will describe our canonical ring ${\mathcal R}$ in
three different
formats introduced by D. Dicks and M. Reid. It turns out that  in
order to deform the
surfaces with $K^2=7$, $p_g=4$ whose canonical system has one base point
and is birational
the third format ($=$  antisymmetric and extrasymmetric format) is the most
suitable one.

The result for this case is

\begin{theo}
Let $S$ be a minimal (smooth, connected) surface with $K^2 = 7$ and
$p_g(S) = 4$, whose
canonical system $|K_S|$ has one (simple) base point
$x \in S$ and yields a birational canonical map.
Then the canonical ring of $S$ can be presented as $\mathcal{R} =
\mathbb{C}[y_0,y_1,y_2,y_3,w_0,w_1,u]/
\mathcal{I}$, where deg $(y_i,w_j,u) = (1,2,3)$
and  the ideal of relations
$\mathcal{I}$ of $\mathcal{R}$ is generated by the
$4 \times 4$ - pfaffians of the following
 antisymmetric and extrasymmetric matrix

$$
P =
\left(
\begin{array}{cccccc}
0 & 0 & w_0 & Q(y_0,y_1,y_2) & w_1 & u \\
  & 0 & y_1 & y_3 & y_2 & w_1 \\
  &  & 0 & -u+\overline{C} & y_3\overline{Q}_3(w_0,y_0,y_1,y_3) & Q\overline{Q}_3 \\
  & &  & 0 & y_1\overline{Q}_1(w_0,w_1,y_0,y_1,y_3) & w_0\overline{Q}_1 \\
& &  &  & 0 & 0 \\
-sym &  &  &  & & 0
\end{array}
\right),
$$

where $Q, \overline{Q}_1, \overline{Q}_3$ are quadratic
forms of a subset of the given variables as indicated,
and $\overline{C}$ is a cubic form. Moreover $\overline{C}$  does not depend on $u$ and
on the  $w_iy_j$'s for
$j
\leq 1$, $\overline{Q}_1$ does not depend on $y_3^2$.
\end{theo}

In the fourth chapter we will finally show how the above
presentation of the canonical ring allows a deformation
to the canonical ring $\mathcal{S}$ of a surface with $K^2 = 7$, $p_g
= 4$ and free canonical
system.

These last surfaces and their canonical rings are described in
\cite{CatBu}, and it follows also from (\cite{BuEi})
that in this presentation the relations can be given by the $4
\times 4$ - Pfaffians of a
skew-symmetric $5 \times 5$ matrix.

Our final  result is

\begin{theo}

Let $P$ be an  antisymmetric and extrasymmetric  matrix as in the
previous theorem. Consider the 1-parameter
family of rings
${\mathcal S}_t={\mathbb C}[y_0,y_1,y_2,y_3,w_0.w_1,u]/{\mathcal
I}_t$ where the ideal
${\mathcal I}_t$ is given by the $4 \times 4$ pfaffians  of the
 antisymmetric and extrasymmetric matrix
$$
P_t =P+
\left(
\begin{array}{cccccc}
0 & t &0 & 0 & 0 & 0 \\
  & 0 & 0 & 0 & 0 & 0 \\
  &  & 0 & 0 & 0 & 0 \\
  &  &  & 0 & 0 & 0 \\
  &  & &  & 0 & t\overline{Q}_1\overline{Q}_3 \\
  -sym &  &  &  & & 0
\end{array}
\right).
$$

This is a flat family and  describes a flat deformation of the
surface corresponding to the matrix $P$ to surfaces with $p_g=4$, $K^2=7$ and
with $|K|$ base point free.
For $t \neq 0$,
${\mathcal S}_t$ is isomorphic to ${\mathbb
  C}[y_0,y_1,y_2,y_3,w_0.w_1]/J_t$
where $J_t$ is the ideal generated by the $4 \times 4$
pfaffians of  the
matrix
$$
\left(
\begin{array}{ccccc}
0 & y_1 & -y_3 & -t^2w_0 & -tQ \\
   & 0 & y_2 & t^3\overline{Q}_3 & tw_1 \\
   &  & 0 & t^2w_1 & t^2\overline{Q}_1  \\
   &  &  & 0 & -t^3c-t^2y_1Q+t^2y_3w_0\\
-sym &  &  &  & 0
\end{array}
\right).
$$
\end{theo}

Then we compare the above family with the one predicted by Enriques
showing (as mentioned above) that it is a completely different type of
degeneration.

Finally we have two appendices. In fact, although theoretically all
the computation can be done by hands, it is better to use a computer
program (as we did with Macaulay 2) to shorten the time needed and be
sure that no computation's mistakes occured. We have put in the
appendix the two Macaulay 2 scripts (without output) that
we needed. The second appendix is interesting, because it shows how
the computer suggested to us the $5 \times 5$ matrix appearing 
in Theorem  0.3.

\section{Canonical systems with base points}

Let $S$ be a minimal surface of general type defined over the complex
numbers and let $|K_S|$ be
its canonical system. If $H^0(S,\mathcal{O}_S(K_S)) \neq 0$, then
$|K_S|$ defines a rational map

$$
\varphi_{|K|} : S -- \rightarrow \mathbb{P}^{p_g-1},
$$

where $p_g = p_g(S) := \dim H^0(S,\mathcal{O}_S(K_S))$ is the geometric
genus of $S$.

\smallskip
Throughout this paper we make the following

\medskip
{\bf Assumption.} $|K_S|$ has no fixed part.

\medskip
Let $\pi : \tilde{S} \longrightarrow S$ be a (minimal) sequence of
blow - ups such that the movable
part $|H|$ of $|\pi^*K_S|$ has no base points. Then we have a
commutative diagram:

\medskip

$$
\begin{array}{ccccl}
\tilde{S} & \multicolumn{3}{c}{\vrpf{30}{\scriptstyle{\varphi_{|H|}}}} &
   \mathbb{P}^{p_g-1} \\
& \rupfl{\scriptstyle{\pi}} & & \ropfr{\scriptstyle{\varphi_{|K_S|}}} & \\
& & S & &  \\
\end{array}
$$

\medskip
Since $\pi$  is a sequence of blow ups $\pi_i : S_i \rightarrow S_{i-1}$
with centre a point
$p_i \in S_{i-1}$, we denote by $E_i$ the $(-1)$-divisor in $\tilde{S}$
given by the full transform of $p_i$, and we denote by $m_i$ the multiplicity
in $p_i$ of the proper transform of a general divisor in $K_S$, so that

\begin{rem}
a) $K_{\tilde{S}} \equiv H + \sum_i \ (m_i +1) E_i$. \\
To simplify the notation, we set $ \mathcal{E} := \sum_i  m_i E_i ,\
E := \sum_i E_i$.\\

b) If $\varphi_{|K_S|}$ is not composed with a pencil, then
$\varphi_{|H|} : \tilde{S}
\longrightarrow \Sigma_1 \subset \mathbb{P}^{p_g-1}$ is a
generically finite morphism
from
$\tilde{S}$ onto a surface $\Sigma_1$ in $\mathbb{P}^{p_g-1}$ of degree
$d = H^2 = K_S^2 - \sum_i \ m_i^2$ and

$$
\varphi_{|H|}^*(\mathcal{O}_{\Sigma_1}(1)) = \mathcal{O}_{\tilde{S}}(H)
=
\mathcal{O}_{\tilde{S}}(\pi^* K_S - \mathcal{E} )
=
\mathcal{O}_{\tilde{S}}(K_{\tilde{S}} - \mathcal{E} - E).
$$

\end{rem}

\begin{defi}
\

a) Let us denote by $\mathcal{F}_1$ the coherent sheaf
of $\mathcal{O}_{\Sigma_1}$ -
modules  $(\varphi_H)_*\mathcal{O}_{\tilde{S}}$. \\

b) We define $\tilde{\mathcal{R}}(S)$ as the graded ring associated to the
divisor $H$ on $\tilde{S}$, thus
$$
  ~\tilde{\mathcal{R}}(S) :=
\bigoplus_{m=0}^{\infty} H^0(\tilde{S}, \mathcal{O}_{\tilde{S}}(mH))=
\bigoplus_{m=0}^{\infty} H^0(\Sigma,
\mathcal{F}_1(m)) .
$$
\end{defi}

We make the following easy observation:

\begin{rem}
$\tilde{\mathcal{R}}(S)$ is a (graded) subring of the canonical ring
$\mathcal{R}(S) :=
\bigoplus_{m=0}^{\infty} H^0(S, \mathcal{O}_S(mK_S))$.
\end{rem}

{\em Proof.} The claim follows from the fact that the natural homomorphism

$$
H^0(\tilde{S}, \mathcal{O}_{\tilde{S}}(mH)) \longrightarrow H^0(\tilde{S},
\mathcal{O}_{\tilde{S}}(mK_{\tilde{S}})) \cong H^0(S, \mathcal{O}_S(mK_S))
$$

is injective for all $m \geq 0$. 

\hfill $\underline{Q.E.D.}$

\

\begin{rem}
We have by our assumption

$$
\tilde{\mathcal{R}}_1 = H^0(\tilde{S}, \mathcal{O}_{\tilde{S}}(H)) =
H^0(\tilde{S},
\mathcal{O}_{\tilde{S}}(K_{\tilde{S}})) = \mathcal{R}_1.
$$

\end{rem}

\medskip
Consider first the Stein factorization of $\varphi_{|H|}$:

\medskip

$$
\begin{array}{ccccl}
\tilde{S} & \multicolumn{3}{c}{\vrpf{30}{\scriptstyle{\varphi_{|H|}}}} &
   \Sigma_1 \subset \mathbb{P}^{p_g - 1} \\
& \rupfl{\scriptstyle{\delta}} & & \ropfr{\scriptstyle{\epsilon_1}} & \\
& & Y  & &  \\
\end{array}
$$

where in general $Y$ is a normal algebraic surface, $\delta$ has connected
fibres and $\epsilon_1$ is  a finite morphism.

\medskip
We shall moreover from now on make the following assumption\\

{\em B) $\varphi_{|H|}$ is a birational morphism onto its image}, whence in
particular $p_g(S) \geq 4$.
\\

Under the above assumption we shall moreover consider a general projection of
$\Sigma_1$ to a surface $\Sigma$ in $\mathbb{P}^3$.

We have therefore the following diagram

$$
\begin{array}{ccccl}
\tilde{S} & \multicolumn{3}{c}{\vrpf{30}{\scriptstyle{\varphi}}} &
   \Sigma \subset \mathbb{P}^{3} \\
& \rupfl{\scriptstyle{\delta}} & & \ropfr{\scriptstyle{\epsilon}} & \\
& & Y  & &  \\
\end{array}
$$

and we may therefore assume

{\em B') $\varphi:\tilde{S} \rightarrow \Sigma$ is a birational morphism}.  \\

We shall write the singular locus
$Sing(\Sigma)$ of
$\Sigma$ as
$\Gamma
\cup Z$, where $\Gamma$ is the subscheme
of $\Sigma$ corresponding to the conductor ideal ${\mathcal C}$
of the normalization morphism $\epsilon$ and where
$Z$ is the finite set  $ \epsilon (Sing(Y)) \subset \Sigma$
(note that if the support of $\Gamma$ is disjoint from $Z$, also $Z$ has a
natural subscheme structure given by the adjunction ideal).

\medskip

$\Sigma$ is Cohen-Macaulay, whence $\Gamma$ is a pure subscheme
of codimension $1$.

\smallskip
We  remark that our methods apply also if the degree of
$\varphi$ is equal to two,
  but we have then to make more complicated technical assumptions.

\medskip

Defining $\mathcal{F}:
= (\varphi)_*\mathcal{O}_{\tilde{S}}$, we have
$$
  ~\tilde{\mathcal{R}}(S) :=
\bigoplus_{m=0}^{\infty} H^0(\tilde{S}, \mathcal{O}_{\tilde{S}}(mH))=
\bigoplus_{m=0}^{\infty} H^0(\Sigma,
\mathcal{F}(m)) ,
$$
whence we may observe that the graded ring $\tilde{\mathcal R}$ is
  a module over the polynomial ring $\mathcal A :=
\mathbb{C}[y_0,y_1,y_2,y_3] $ (the homogeneous coordinate ring of
$\mathbb{P}^3$).

Since this module has a support of codimension $1$, it has a graded free
resolution of length equal to $1$ if and only if it is a Cohen-Macaulay
module.

The following result is essentially the same result as theorem 2.5 of
\cite{Cil1}.
\begin{prop}\label{CMpn}
$\tilde{\mathcal R}$ is
  a Cohen Macaulay $\mathcal A $-module if and only if the subscheme
  $\Gamma \subset \mathbb{P}^3$ is projectively normal.
\end{prop}

{\it Proof.}
It is well known that the Cohen Macaulay property is equivalent to the
vanishing of the cohomology groups

$$
H^1(\Sigma,\mathcal{F}(m)) = H^1(\Sigma, \epsilon_* \mathcal{O}_Y
\otimes \mathcal{O}_{\Sigma}(m)) =
H^1(Y, \mathcal{O}_Y(mH)),
$$
for all $m$.

By Serre's theorem B($m$) these groups obviously vanish for $m >>0$.

Serre-Grothendieck duality tells us that these are the dual vector spaces of
$ Ext^1(\epsilon_* \mathcal{O}_Y(m),\omega_{\Sigma})$.

By the local-to-global spectral sequence of the Ext groups, there is an exact
sequence
\begin{multline*}
0
\rightarrow
H^1({\mathcal H}om(\epsilon_* \mathcal{O}_Y(m),\omega_{\Sigma}))
\rightarrow
Ext^1(\epsilon_* \mathcal{O}_Y(m),\omega_{\Sigma})
\rightarrow\\
H^0({\mathcal E}xt^1(\epsilon_* \mathcal{O}_Y(m),\omega_{\Sigma}))
\rightarrow
H^2({\mathcal H}om(\epsilon_* \mathcal{O}_Y(m),\omega_{\Sigma}))
.
\end{multline*}

But ${\mathcal E}xt^1(\epsilon_* \mathcal{O}_Y(m),\omega_{\Sigma})$ is
zero because $\epsilon_* \mathcal{O}_Y$ is a Cohen-Macaulay
$\mathcal{O}_{\Sigma}-$module and $ \omega_{\Sigma}$ is invertible.

Therefore  it follows that $Ext^1(\epsilon_*
\mathcal{O}_Y(m),\omega_{\Sigma})=0$ if and only if  $H^1({\mathcal
H}om(\epsilon_*
\mathcal{O}_Y(m),\omega_{\Sigma}))=0$.

In turn, since $\omega_{\Sigma}={\mathcal O}_{\Sigma}(d-4)$,
$H^1({\mathcal H}om(\epsilon_* \mathcal{O}_Y(m),\omega_{\Sigma}))=
H^1({\mathcal C}(m+d-4))$, where ${\mathcal C}$ is the conductor ideal of
$\epsilon$, and this last group, in view of the exact sequence
$$
0
\rightarrow
{\mathcal C}
\rightarrow
{\mathcal O}_{\Sigma}
\rightarrow
{\mathcal O}_{\Gamma}
\rightarrow
0
$$
is the cokernel of the map $H^0({\mathcal O}_{\Sigma}(m+d-4))
\rightarrow
H^0({\mathcal O}_{\Gamma}(m+d-4))
$.

Since however $H^0({\mathcal O}_{\Sigma}(n))$ is a quotient  of
$H^0({\mathcal O}_{{\mathbb P}^3}(n))$ we conclude that our desired vanishing
is equivalent to the projective normality of $\Gamma$.

\hfill $\underline{Q.E.D.}$

\

Recalling that
$ \mathcal{O}_{\tilde{S}}(H)
=
\mathcal{O}_{\tilde{S}}(\pi^* K_S - \mathcal{E} )$,
  we consider now the exact sequence

$$
0 \longrightarrow \mathcal{O}_{\tilde{S}}(mH) \longrightarrow
\mathcal{O}_{\tilde{S}}(m \pi^* K_S ) \longrightarrow
\mathcal{O}_{m\mathcal{E}} \longrightarrow 0.
$$

Observe moreover that, since  $S$ is minimal and of general type,
for $m \geq 2$ we have

$$
H^1(\tilde{S}, \mathcal{O}_{\tilde{S}}(m \pi^* K_S) =
H^1(S, \mathcal{O}_S(mK_S)) = 0.
$$

Whence,
we arrive for each $m \geq 2$ to
the following crucial exact
sequence, which will be repeatedly used  in the sequel.

$$
0 \longrightarrow \tilde{\mathcal{R}}_m \longrightarrow \mathcal{R}_m
\longrightarrow
H^0(\tilde{S}, \mathcal{O}_{m\mathcal{E}}) \longrightarrow H^1(\tilde{S},
\mathcal{O}_{\tilde{S}}(mH))  \longrightarrow 0. \leqno(ii)
$$

The vanishing $H^1(\Sigma, \mathcal{F}(m)) = 0$ for all
$m$, is clearly equivalent to the chain of equalities:

$$
dimH^1(\tilde{S}, \mathcal{O}_{\tilde{S}}(mH)) = dimH^0(\Sigma,
R^1(\varphi_H)_*
\mathcal{O}_{\tilde{S}}(m)) = length (R^1(\varphi_H)_*
\mathcal{O}_{\tilde{S}}) := l.
$$

\medskip
Putting together the above considerations we obtain the following

\medskip
\begin{rem}\label{dimrtilda}
$\tilde{\mathcal{R}} $ is a Cohen Macaulay module over the polynomial ring
  $\mathcal{A}$ if and only if the surface $S$ is regular ($H^1(S,
\mathcal{O}_S(K_S)) = 0$) and

$ dim {\mathcal{R}}_m - dim
\tilde{\mathcal{R}}_m = dim H^0(\tilde{S}, \mathcal{O}_{m\mathcal{E}}) - l =
\sum_i \frac{mm_i(mm_i+1)}{2} -l$.
\end{rem}

{\it Proof.}

Assume $\tilde{\mathcal R}$  to be Cohen Macaulay:

since we know that $H^1({\mathcal O}_{\mathcal{E}})=0$, and that the map
$\tilde{\mathcal R}_1
\rightarrow {\mathcal R}_1$ is an isomorphism (from the definition of
$\tilde{\mathcal R}$),  we get an exact sequence
$$
0
\longrightarrow
H^0(\tilde{S}, \mathcal{O}_{\mathcal{E}}) \longrightarrow H^1(\tilde{S},
\mathcal{O}_{\tilde{S}}(H))  \longrightarrow H^1(S,
\mathcal{O}_{S}(K_S)) \longrightarrow 0.
$$
which shows that the surface must be regular.

If conversely the surface is regular,
$H^1(S,\mathcal{O}_S)=H^1(S,\mathcal{O}_S(K_S))=0$ , so the sequence $ii)$
is exact also for $m=0,1$, then $\forall m \in {\mathbb Z}$.

In particular, $ dim {\mathcal{R}}_m - dim
\tilde{\mathcal{R}}_m = dim H^0(\tilde{S}, \mathcal{O}_{m\mathcal{E}}) -
dim H^1(\tilde{S}, \mathcal{O}_{\tilde{S}}(mH))$, whence
$\tilde{\mathcal R}$ is a Cohen Macaulay
${\mathcal A}-$module iff  $H^1(\Sigma, \mathcal{F}(m)) = 0$ , i.e., iff
$dim H^1(\tilde{S}, \mathcal{O}_{\tilde{S}}(mH))=l$, equivalently,  iff
$ dim {\mathcal{R}}_m - dim
\tilde{\mathcal{R}}_m = dim H^0(\tilde{S}, \mathcal{O}_{m\mathcal{E}}) - l =
\sum_i \frac{mm_i(mm_i+1)}{2} -l$. 
 
\hfill $\underline{Q.E.D.}$

\

\begin{rem}
$\tilde{\mathcal{R}} \subset \mathcal{R}$ is a subring, but it is
easy to see that $\mathcal{R}$
is not a finitely generated $\tilde{\mathcal{R}}$ - module.
\end{rem}

We assume now that $\tilde{\mathcal R}$ is a Cohen-Macaulay ${\mathcal
A}$-module, and we observe that it contains the coordinate ring of
$\Sigma$ and is
contained in ${\mathcal R}$.

Therefore, choosing a
minimal system of generators
$v_1=1,v_2,
\ldots v_n$ of $\tilde{\mathcal R}$ as an ${\mathcal A}-$module, defining
$l_i:= \deg v_i$ we find (by Hilbert's syzygy theorem, as in \cite{CatBu}), a
resolution of the form
$$
(\#)\ \  \ \  \  \ \ \ \
0
\longrightarrow
\oplus_{j=1}^h {\mathcal A}(-r_j)
\stackrel{\alpha}{\longrightarrow}
\oplus_{i=1}^h {\mathcal A}(-l_i)
\longrightarrow
\tilde{\mathcal R}
\longrightarrow
0;
$$
Under the assumption that $\varphi$ be birational follows that $\Sigma$
has equation $f := \det \alpha=0$.

As in \cite{CatBu}, being $\tilde{\mathcal{R}}$ a ring, the matrix has
to fulfill the standard Rank Condition, which we will later recall.

In order to describe the ring ${\mathcal R}$, we first look for
generators of ${\mathcal R}$ as an ${\mathcal A}$-module.

On the other hand, when looking for generators of ${\mathcal R}$ as a ring,
we may restrict ourselves to consider elements of low degree by virtue of the
following result by M. Reid (cf. \cite{Rei}, cf. also Ciliberto
\cite{ciliberto}).

\begin{theo}\label{boundongens}
Let $X$ be a canonical surface (i.e., the canonical model of a
surface of general type). We
suppose that \\
(i) $p_g(X) \geq 2$, $K_X^2 \geq 3$, \\
(ii) $q(X) = 0$, \\
(iii) $X$ has an irreducible canonical curve $C \in |K_X|$. \\
Then the canonical ring $\mathcal{R} = \mathcal{R}(X, K_X)$ of $X$ is
generated in degrees $\leq 3$
and its relations are generated in degrees $\leq 6$.
\end{theo}

Now we define the $\tilde{\mathcal{R}}$-module
$$
M: = \Gamma_*({\mathcal C} \omega_{\Sigma}) =\Gamma_* ( \omega_Y),
$$

where ${\mathcal C}$ is the conductor ideal of $\epsilon$ and
$\Gamma_*({\mathcal F})$ denotes as usual $\oplus_{n \in {\mathbb Z}}
H^0({\mathcal F}(n))$.
We consider the following chain of inclusions of ${\mathcal A}-$modules

$$
{\mathcal A}/(f)
\subset
\tilde{\mathcal R}
\subset
\Gamma_*(\varphi_* \omega_{\tilde{S}})[-1]
\subset
M[-1].
$$

We observed that $\tilde{\mathcal R}$ is Cohen Macaulay if and only if it
has a free resolution as an ${\mathcal A}-$module of the form $(\#)$

$$
0
\rightarrow
L_1
\stackrel{\alpha}{\rightarrow}
L_0
\rightarrow
\tilde{\mathcal R}
\rightarrow
0.
$$

Dualizing it, we get

$$
0
\rightarrow
L_0^{\vee} \otimes \Gamma_*(\omega_{{\mathbb P}^3})
\stackrel{\alpha^t}{\rightarrow}
L_1^{\vee} \otimes \Gamma_*(\omega_{{\mathbb P}^3})
\rightarrow
Ext^1_{\mathcal A} (\tilde{\mathcal R}, \Gamma_*(\omega_{{\mathbb P}^3}))
\rightarrow
0.
$$

By virtue of the exact sequence
$$
0
\rightarrow
\Gamma_*(\omega_{{\mathbb P}^3})
\stackrel{f}{\rightarrow}
\Gamma_*(\omega_{{\mathbb P}^3})(\deg f)
\rightarrow
\Gamma_* ( \omega_{\Sigma})
\rightarrow
0,
$$
and since $Ext^1_{\mathcal A} (\tilde{\mathcal R}, f)=0$, we get
$$Ext^1_{\mathcal A} (\tilde{\mathcal R}, \Gamma_*(\omega_{{\mathbb P}^3}))=
Hom_{\mathcal A} (\tilde{\mathcal R}, \Gamma_*(\omega_{\Sigma}))=
\Gamma_*({\mathcal H}om_{{\mathcal O}_{\Sigma}} ({\mathcal F},
\omega_{\Sigma}))=
\Gamma_*({\mathcal C} \omega_{\Sigma})=M.
$$

Moreover, $M$ satisfies the Ring Condition (cf. \cite{djvs})

$$
\tilde{\mathcal R}= Hom(M, {\mathcal A}/(f))=Hom (M,M)
$$
or, in other words, there is a bilinear pairing $\tilde{\mathcal R} \times M
\rightarrow
M$ which, in the given bases, is determined by the matrix $\beta=
\Lambda^{h-1}(\alpha)$. In turn the Ring Condition is equivalent to the so
called Rank Condition  for $(\alpha)$ :
there exist elements
$\lambda_{jh}^k$ of $\mathcal A$ such that
$\beta_{1k}=\sum \lambda_{jh}^k \beta_{jh}$.

Since we have the inclusion $\tilde{\mathcal R} \subset M[-1]$, we can fix
bases
$v_1=1,v_2, \ldots, v_h$, for
$\tilde{\mathcal R}$,
$z_1,z_2, \ldots, z_h$, for $M$, with
$v_iz_j=\frac{\beta_{ij}z_1}{\beta_{11}}$ and such that $z_1$ is the image of
$v_1=1$.

For the same reason, $v_2, \ldots, v_n$ can be written as linear
combinations $v_i = \sum \zeta_{ij} z_j/z_1$.

Now, as in \cite{CatBu}, the ring structure of $\tilde{\mathcal R}$ is
equivalent to the Rank Condition  which can also be phrased as follows:

\

{\bf R.C.} the ideal of the $(h-1) \times (h-1)$ minors of $\alpha$
coincide with the ideal of the $(h-1) \times (h-1)$ minors of
the matrix $\alpha'$ obtained by  deleting the first
row of $\alpha$.

\

The polynomials $\lambda_{jh}^k$ with
$\beta_{1k}=\sum \lambda_{jh}^k \beta_{jh}$ determine therefore the ring
structure of $\tilde{\mathcal R}$ by the following multiplication rule :\\

$$v_iv_h=\sum_{j,k}
\zeta_{ij}
\lambda_{jh}^k v_k.$$

We have now all the general ingredients at our disposal and we can explicitly
describe our method in order to compute the canonical ring of the regular
surfaces of general type with given values of the invariants $K^2$, $p_g$,
canonical system with base points (and without fixed part),
and birational canonical map.

Under the above assumptions, observe that obviously $p_g \geq 4$, moreover
  Castelnuovo's inequality $K^2 \geq 3p_g-7$ holds, in particular the
hypotheses of theorem
\ref{boundongens} are  fullfilled.

We need  first of all to assume that the subscheme $\Gamma$ of
${\mathbb P}^3$ given by the conductor ideal of the normalization of $\Sigma$
be projectively normal.

The last assumption, as we just saw, ensures that the ring $\tilde{\mathcal
R}$ is a Cohen-Macaulay ${\mathcal A}-$module; whence, argueing as in
\cite{CatBu} we can find a length $1$ presentation of $\tilde{\mathcal
  R}$ as an ${\mathcal A}-$module, given by a square matrix $\alpha$
fullfilling the Rank Condition .

Let $v_1, \ldots, v_h$ be the generators of $\tilde{\mathcal R}$ we
used in order to write down $\alpha$, and let $z_1, \ldots z_h$
be the dual generators of $M$ (i.e., the module $M$ is generated by the $z_i$'s
and presented by the matrix $\alpha^t$).

We have seen that there is an inclusion $\tilde{\mathcal R} \subset
 M[-1]$; assuming by sake of simplicity that $\omega_Y$ is Cartier,
 one can write explicit sections $\sigma_d$ ($d \in {\mathbb N}$) of
 suitable line bundles ${\mathcal L}_d$ so that the above inclusion is
 obtained multiplying every element $r_d$ homogeneous of degree $d$ in
 $\tilde{\mathcal R}$, by
 $\sigma_d$. The $\sigma_d$'s are of the form $e_d \cdot c$ where $c$ is
 a section of the dual of the relative canonical bundle of the map
 $\delta : \tilde{S} \rightarrow Y$, and $e_d$ is supported on the
 exceptional locus of $\pi:\tilde{S} \rightarrow S$.

It is not possible to construct a similar inclusion ${\mathcal
  R}\subset M$, but we can consider the module $\oplus_n
  H^0(K_{\tilde{S}}+nH)$. This is the submodule of $M$ given by the
  elements divisible by $c$ (in particular it contains $\tilde{R}$),
  and clearly it is a submodule of ${\mathcal R}$ so it is completely
  natural to denote it by 
$M \cap {\mathcal R}$. 

The second step of our method is to study
  this module: first we compute the subset $\{w_1, \ldots w_r\}$ of a set of
  generators for $M \cap {\mathcal R}$ as an ${\mathcal A}$-module,
  consisting of the elements of  degree $\leq 3$ (with the
  grading of ${\mathcal R}$). Then wefind the relations holding among
  them in degree $\leq 6$. In fact, by theorem \ref{boundongens}, generators
  and relations in higher degrees will not   be relevant for the
  canonical ring. 

The elements $y_i's$, $w_j$'s will not in general generate the
canonical ring; the third step of our method consists in the
research of the missing generators and relations.  
 
It is clear that in every case ${\mathcal R}_1 \subset
{\mathcal R} \cap M$. We assume now that the base points are simple
(but a similar analysis can be carried out in every case), i.e.
$$
{\mathcal E}=E. \leqno{(*)}
$$

Directly by the definition follows the equality
$K_{\tilde S}+H= 2 \pi^* K_S$, so
$H^0(K_{\tilde S}+H)=H^0(2K_S)$, i.e. ${\mathcal R}_2 \subset
{\mathcal R} \cap M$. Instead, if ${\mathcal E} \neq
E$ the equality $h^0(K_{\tilde S}+H)=h^0(2K_S)$ cannot hold, otherwise
$|2K_S|$ would not be base point free, as it has to be for $p_g >0$
(\cite{Francia}).

In degree $3$ we have the exact sequence
$$
0
\rightarrow
H^0(K_{\tilde{S}}+2H)
\rightarrow
H^0(3K_{S})
\rightarrow
H^0({\mathcal O}_{\mathcal E})
\rightarrow
0
$$
since $H^1(K_{\tilde{S}}+2H)=0$ (by Mumford's vanishing theorem).

Therefore, ${\mathcal R}_3 \cap M$ has codimension $3$ in ${\mathcal
  R}_3$, i.e.   
we need $l$ elements $u_1, \ldots, u_l$, to complete a
basis of $H^0(K_{\tilde{S}}+2H)$ to a basis of $H^0(3K_{S})$.

By theorem $\ref{boundongens}$ the generators $y_0,\ldots ,y_3$ of
${\mathcal A}$ (seen as elements of $H^0(K_{\tilde S})$, 
together with $w_1, \ldots, w_r$ and $ u_1, \ldots, u_l$ 
are a system of generators of ${\mathcal R}$ as a ring.

The relations as ${\mathcal A}$-module among the $v_i$'s and among
the $w_j$'s are determined by the matrix $\alpha$, and similarly the
relations  given by the products of type $v_iv_j$, and $v_iw_j$ are
also determined by $\alpha$: these provide automatically 
a list of relations among the above generators of ${\mathcal R}$.

Some relation is still missing; in particular our method do not
produce automatically any relation involving the $u_k$'s. In order to
complete the analysis one needs to find a way to espress the $u_k$'s
``in terms of $M$'', so that the known relation among the elements
of ${\tilde R}$ and $M$ will produce also relations involving
them. This should be possible case by case by ``ad hoc'' arguments,
but we do not know a general argument.

We devote the next sections to the
application of the above method to finding a description of the canonical
ring of the stratum
of the moduli space of surfaces of general type with
$K^2=7$ and $p_g=4$ corresponding to surfaces with
birational canonical map and whose canonical system has base points.

\section{The canonical ring of surfaces with $K^2 = 7$, $p_g =4$
birational to a sextic:
from  algebra to geometry}

Let $S$ be a minimal (smooth, connected) surface with $K^2 = 7$ and
$p_g(S) = 4$. \\
We remark that $S$ is automatically regular (cf. \cite{O. Debarre}),
i.e. $q(S) = 0$. \\

In \cite{Ba} the first author gave an exact description of minimal
surfaces with $K^2 = 7$ and
$p_g(S) = 4$, where the canonical system has  base points, proving
in particular that if moreover the canonical map is
birational
then the canonical system $|K_S|$ has exactly one simple base point
$x \in S$.

Let $\pi : \tilde{S} \longrightarrow S$ be the blow up of $S$ in $x$
and let $E := \pi^{-1}(x)$ be the exceptional curve of $\pi$. Thus we
have:

$$
|K_{\tilde{S}}| = |\pi^*K_S| + E = |H| + 2E,
$$

where $|H|$ is base point free.\\
Thus we will assume in this paragraph that $\varphi_{|H|} : \tilde{S}
\longrightarrow \Sigma := \{F_6 = 0 \}$ is a birational morphism
(from $\tilde{S}$ onto a surface $\Sigma$ of  degree six in $\mathbb{P}^3$).

\medskip

We recall now the  description of the
minimal surfaces with $K^2 = 7$ and
$p_g(S) = 4$, whose  canonical system has exactly one base point
and whose canonical map is
birational given in \cite{Ba}.\\

First we need the following definition.

\begin{defi}
A {\em generalized tacnode} is a two dimensional elliptic
hypersurface singularity
$(X,0)$, such that the fundamental cycle has self intersection $(-2)$. \\
In particular, $(X,0)$ is Gorenstein and by \cite{Lau}, theorem
$(1.3)$, $(X,0)$ is a
double point singularity, whose  local analytic equation is given by

$$
z^2 = g(x,y),
$$

where $g$ vanishes of order four in $0$. The normal cone of the
singularity is given by the plane $\{z=0\}$, called the tacnodal plane.
\end{defi}

More precisely, a generalized tacnode $(X,0)$ is the singularity obtained as
the double cover branched along a
curve with a quadruple point, which after a blow up decomposes in at
most simple triple points
or double points.

\begin{theo}\label{geometria} (\cite{Ba}, theorem $3.11$, $5.5$)
1) Let $S$ be a minimal surface with $K^2 = 7$ and
$p_g(S) = 4$, whose  canonical system has exactly one base point
and whose canonical map is
birational. Then the blow-up $\tilde{S}$ of the base point is the minimal
desingularization of a surface $\Sigma \subset \mathbb{P}^3$ of degree six
with the following properties: \\
(a) ~the double curve $\Gamma \subset \Sigma$ is a plane conic, \\
(b) ~if $\gamma \subset \mathbb{P}^3$ is the plane containing
$\Gamma$, then $\Sigma$
has a generalized tacnode $o \in \gamma \backslash \Gamma$ with tacnodal plane
$\alpha \neq \gamma$, \\
(c)  the image $\varphi_H(E) $ of the exceptional curve equals
the line $\alpha \cap
\gamma$.
\\ 2) The surfaces with $K^2 = 7$ and $p_g = 4$ such that the canonical
system has exactly one base point and $\varphi_{|K|}$ is birational form an
irreducible set $\mathfrak{M}_{(I.1)}$ of dimension $35$ in their moduli space.
\end{theo}

Moreover, it was shown (ibidem) that for a general element of
$\mathfrak{M}_{(I.1)}$ the canonical image
$\Sigma$ has an equation of the form

$$
\alpha^2Q^2 + \gamma F_5 = 0,
$$

where $F_5$ is an element of the linear subsystem $\Delta \subset |5H
- \Gamma - 2o|$ in
$\mathbb{P}^3$ consisting of  quintics with tangent cone $\alpha^2$ in $o$ and
where $Q
\subset \mathbb{P}^3$ is an
irreducible quadric containing $\Gamma$.

In this section we will study the ring $\tilde{\mathcal R}$ and the module
$M$ defined in the previous section in the case of the surfaces in the
above class.

This study will allow us to compute the canonical ring and to give
a purely algebraic proof of part 1) of theorem \ref{geometria},
under a few generality assumptions.

For the convenience of the reader we will give here a list of
notation (partly already
introduced in the last section) which will be frequently used in the
following.

\smallskip
{\bf Notation:} \\
We consider:
\begin{itemize}
\item $\pi : \tilde{S} \rightarrow S$, the blow-up of the base point
of $|K_S|$;

\item $\varphi : \tilde{S} \rightarrow
\! \! \! \! \!
\rightarrow \Sigma \subset {\mathbb P}^3$, the morphism induced by
the canonical system;

\item $\epsilon : Y \rightarrow
  \Sigma$, the normalization of $\Sigma$;

\item $\delta : \tilde{S} \rightarrow Y$, such that $\varphi=\epsilon
  \circ \delta$;

\item $e$, a generator of $H^0(\tilde{S}, \mathcal{O}_{\tilde{S}}(K_{\tilde{S}} -
\pi^*K_S))$ and $E$
the corresponding divisor;

\item $c$, a generator of $H^0(\tilde{S}, \mathcal{O}_{\tilde{S}}(\delta^* K_Y -
K_{\tilde{S}}))$ and $C$ the
corresponding divisor (this makes sense under the assumption 2) below);

\item $F_6$ an equation of $\Sigma \subset \mathbb{P}^3$.

\smallskip
We have the following list of graded rings respectively modules:

\item $
\mathcal{A} := \bigoplus_{m=0}^{\infty}
S^m(H^0(\tilde{S},\mathcal{O}_{\tilde{S}}(H))) =
{\mathbb C}[x_0,x_1,x_2,x_3],
$

\item $
\tilde{\mathcal{R}} := \bigoplus_{m=0}^{\infty}
H^0(\Sigma,(\varphi_{|H|})_* \mathcal{O}_{\tilde{S}}(m)) =
\bigoplus_{m=0}^{\infty} H^0(\tilde{S}, \mathcal{O}_{\tilde{S}}(mH)),
$

\item $
M := \bigoplus_{m \in {\mathbb Z}}
H^0(\Sigma,{\mathcal C}\omega_{\Sigma}(m)) =
\bigoplus_{m=-1}^{\infty} H^0(\tilde{S}, \mathcal{O}_{\tilde{S}}((m+1)H+2E+C))
$, where

\item $\mathcal{C} :=
\mathcal{H}om_{\mathcal{O}_{\Sigma}}(\epsilon_*\mathcal{O}_Y,\mathcal{O}_{\Sigma
})$
is
the conductor ideal;

\item $
{\mathcal R}=\bigoplus_{m=0}^{\infty}
H^0(S,mK_S) =
\bigoplus_{m=0}^{\infty} H^0(\tilde{S},\mathcal{O}_{\tilde{S}}(m(H+E))
$, the canonical ring of $S$.
\end{itemize}
\bigskip

The assumptions we will make are the following:

1) The conductor ideal of $\epsilon$ defines a projectively normal
   subscheme of ${\mathbb P}^3$;

2) The singular points of $Y$ do not lie in the preimage of the non
    normal locus of $\Sigma$; in particular it follows that $\omega_Y$
    is Cartier.

These two assumptions give in fact no restriction, as it can be shown
with geometrical arguments (\cite{Ba}, remark 3.1. and
prop. 3.6.(v)).

First, by the results of the previous section, we give a presentation
of the ring $\tilde{\mathcal{R}}$ as ${\mathcal A}$-module.

\begin{theo}\label{FisCM}
$\tilde{\mathcal{R}}$ is a Cohen - Macaulay $\mathcal{A}-$module and
has a resolution (as an $\mathcal{A}$-module) as follows

$$
\begin{array}{ccccccccc}
& & \mathcal{A}(-5) & \alpha& \mathcal{A} & & & & \\
0 & \to & \bigoplus & \to & \bigoplus & \to & \tilde{\mathcal{R}} &
\to & 0 . \\
& & \mathcal{A}(-4) & & \mathcal{A}(-3) & & & &
\end{array}
$$
\end{theo}

{\em Proof.} (cf. \cite{CatBu} for similar computations).

By remark \ref{dimrtilda},
dim $\tilde{\mathcal R}_m = \dim {\mathcal R}_m-\frac{m(m+1)}{2}+1$, therefore
by Riemann-Roch's Theorem, for $m \geq 2$, it equals $\chi({\mathcal
   O}_S)+\frac72m(m-1)-\frac{m(m+1)}{2}+1=3m^2-4m+6$.

$\tilde{\mathcal R}$ is Cohen Macaulay by theorem \ref{CMpn} (and
assumption 1)), whence it has a resolution as ${\mathcal A}-$module
of length 1.

By definition the first generator of $\tilde{\mathcal R}$ has degree $0$
($\tilde{\mathcal R}_0 = H^0(\mathcal{O}_{\tilde{S}})$) and will be
denoted by $1$.

The above dimension formula gives us immediately that there are no
other generators in degrees $\leq 2$ and that one more generator (denoted
by $v$) is needed in degree $3$. Moreover, the relations live in
degrees $\geq 4$.

Again by the above dimension formula we get at least one relation in
degree $4$. If there were two independent relations in degree
$4$ they would have the form $x_0 \cdot v=f_4(x_i)\cdot 1$; $x_1\cdot
v=g_4(x_i)\cdot 1$ and this would force a non trivial relation of the form
$(x_1f_4-x_0g_4)\cdot 1=0$ of degree $5$, contradicting that obviously for any $f \in \mathcal{A}$ the
equality $f \cdot 1=0$
implies that $f$ is a multiple of $F_6$.
By the dimension formula therefore there are no new generators in degree 4.

Again counting the dimensions we get a relation in degree $5$; a
straightforward
computation shows that
for all $m$, dim
$\tilde{\mathcal R}_m=$ dim ${\mathcal A}_m$+ dim ${\mathcal A}_{m-3}
-$dim ${\mathcal A}_{m-4} -$dim ${\mathcal
   A}_{m-5}$;  this shows that the resolution has the form
$$
\begin{array}{ccccccccc}
& & \mathcal{A}(-5) & \alpha & \mathcal{A} & & & & \\
0 & \to & \bigoplus & \to & \bigoplus & \to & \tilde{\mathcal{R}} &
\to & 0 . \\
& & \mathcal{A} (-4) & & \mathcal{A}(-3) & & & &\\
& & \bigoplus & & \bigoplus & & & &\\
& & \mathcal{L} & & \mathcal{L} & & & &
\end{array}
$$
where ${\mathcal L}$ is a free module $\oplus \mathcal{A}(-s_i)$, with
$s_i \geq 5$ for all $i$.

The minimality of the resolution ensures that ${\mathcal L}=0$
(there are no non zero constants as coefficients in $\alpha$;
considering the row of $\alpha$ corresponding to the new generator
of maximal degree ($\geq 5$) we get a row of zeroes, contradicting the
injectivity of $\alpha$).

\hfill $\underline{Q.E.D.}$

\

\begin{cor}\label{resM}
$M$ has a resolution as ${\mathcal A}$-module of the form

$$
\begin{array}{ccccccccc}
& & \mathcal{A}(-4) & ^t\alpha& \mathcal{A}(1) & & & & \\
0 & \to & \bigoplus & \to & \bigoplus & \to & M &
\to & 0 . \\
& & \mathcal{A}(-1) & & \mathcal{A} & & & &
\end{array}
$$
\end{cor}

We denote by $z_{-1}$, $z_0$ the generators of
$M$ in the respective degrees $-1$ and $0$;

\begin{rem}\label{z-1}
Notice that $z_{-1}=e^2c$.
\end{rem}

\begin{cor}\label{rc}
Up to a suitable choice of the generators of the ${\mathcal
A}-$modules $\tilde{\mathcal R}, M$, we can write
$$\alpha=
\left(
\begin{array}{cc}
QG+\gamma B & Qq \\
Q & \gamma
\end{array}
\right);
$$
where deg $(\gamma,Q,q,G,B) = (1,2,2,3,4)$, moreover
$$
v z_{-1}=q z_0
$$
$$
v z_0=B z_{-1} -G z_0
$$
\end{cor}

{\it Proof.} By the above resolution of $M$ we know the degrees of
the entries of
$\alpha$. The Rank Condition for $\alpha$ means that the elements of
the first row are in
the ideal generated by the elements of the second row.
To obtain now the desired form of $\alpha$ it suffices to add a
suitable multiple of the
second row to the first one.

The two relations can be easily obtained writing explicitly the
pairing $\tilde{\mathcal R} \times M \rightarrow M$.

\hfill $\underline{Q.E.D.}$

\

In  order to understand the structure of the canonical ring, we have
now to investigate which elements of the module $M$ can be divided by
$c$. In fact, the graded parts of our rings are
related by the following (commutative) diagram:

$$ (*) \ \ \ \ \ \ \
\begin{array}{ccccc}
\tilde{R}_n    && \stackrel{e^2}{\longrightarrow} &&
H^0(\tilde{S}, \mathcal{O}(nH+2E)) \\
\scriptstyle{e^n} \downarrow   && \swarrow && \downarrow \scriptstyle{c} \\
R_n &&&& M_{n-1}
\end{array}
$$
where the diagonal arrow is multication by $e^{n-2}$ ($n \geq 2$).

In order to write down explicitely the ring, we will fix
a basis $x_0,x_1,x_2,x_3$ for ${\mathcal A}_1=\tilde{{\mathcal R}}_1$.
We will denote by $y_i:=ex_i$ the induced elements in ${\mathcal R}$.

\begin{rem} \label{deg2}
First, since $\gamma \neq 0$ (or
$\Sigma=\{ \gamma F_5 = qQ^2 \}$ would be reducible), we set
$x_3:=\gamma$. Moreover we shall assume from now on that
$G=G(x_0,x_1,x_2)$ (this is clearly possible without loss of
generality by corollary \ref{rc}). 

By the resolution of $M$, $\dim_{\mathbb C} M_1=13$ and a basis of $M_1$ is
provided by the 10 elements of the form $x_ix_jz_{-1}$ and
$x_0z_0,x_1z_0,x_2z_0$ (the only relation in
degree $1$ has the form $x_3z_0 =
-Qz_{-1}$).

By Riemann-Roch, $P_2=12$; recalling that $c|z_{-1}$
diagram $(*)$ shows that we can fix the basis
$x_0,x_1,x_2,x_3=\gamma$ of ${\mathcal A}_1$ such that $c|
x_iz_0$ for $i \geq 1$ but $c$ does not divide
$x_0z_0$, and setting
$w_0:=-\frac{x_1z_0}{c}$, $w_1:=-\frac{x_2z_0}{c}$,
$w_0,w_1,y_iy_j$ is a basis of ${\mathcal R}_2$.
\end{rem}

\begin{lemma}\label{deg3}
$R={\mathbb C}[y_0,y_1,y_2,y_3,w_0.w_1,u]/I$ where deg $(y_i,w_j,u)$ =
$(1,2,3)$, and the generators in degree $\leq 3$ of the ideal $I$ are given
by the vanishing of the $2 \times 2$ minors of the following matrix:

$$\left(
\begin{array}{ccc}
y_1 & y_2 & y_3 \\
w_0 & w_1 & Q(y_i)
\end{array}
\right)
$$
where $Q(y_i)$ is obtained from $Q$ replacing every $x_i$ by the
corresponding $y_i$.

\end{lemma}

{\it Proof.}

In the previous remark, we saw that in degree $\leq 2$, ${\mathcal R}$
is generated by the $y_i$'s and  the $w_j$'s.

Every relation in degree $\leq 3$ of ${\mathcal R}$ is of the form
$w_0l_0(y_i)+w_1l_1(y_i)=g_3(y_i)$, where $l_0$ and $l_1$ are linear
forms, $g_3$ is a cubic form. 

Thus we get from it the following
relation in degree $2$ for $M$ as ${\mathcal
   A}-$module: $z_0(x_1l_0(x_i)+x_2l_1(x_i))=z_{-1}g_3(x_i)$. 

By corollaries \ref{resM} and \ref{rc}, must be a 
multiple of the relation $Qz_{-1}+x_3z_0=0$. Thus we immediately see
that this relation is a linear combination of the three 
given by the $2$ by $2$ minors of the matrix in the statement.

We obtain therefore that the subspace of $R_3$ generated by the
monomials in the $y_i$, $w_i$ has dimension $25$. Since $P_3=26$ we
need to add a new generator $u$.

Finally, there are no more generators in degree $\geq 4$ by theorem
\ref{boundongens}.

\hfill $\underline{Q.E.D.}$

\

We noticed in fact at the end of the previous section that if $E={\mathcal
E}$  we need exactly $l=h^0({\mathcal O}_{\mathcal
E})$ new generators of $\mathcal{R}$ in degree $3$ which do not come
from elements in
$M$.  This gives in our case ($l=1$) exactly one new generator in degree
$3$ (the only one
not vanishing in $E$), that is, our generator $u$.

\begin{rem}\label{cnondividex_0^2z_0}
For later use we observe that $c$ does not divide $x_0^2 z_0$.

This holds since otherwise $\frac{x_0^2z_0}{c}$ would yield
(cf. diagram (*)) an element in $H^0(3H+2E)$, whence
$\frac{ex_0^2z_0}{c}$ would be an element of $R_3$ linearly dependent
upon the monomials in the $y_i$'s, $w_i$'s. This however contradicts 
corollary \ref{resM}.
\end{rem}

As pointed out at the end of the previous section, the main problem in
computing the canonical ring is given by the "additional" generators
in degree $3$ (in our case there is only one, namely $u$).

In the next lemma (part 5) we manage to express $u$ "in terms of $M$"; this
will allow us to compute the canonical ring.

\begin{lemma} \label{geo}
Choosing suitable coordinates in $\mathbb{P}^3$ and suitable
generators of $M$ and
$\tilde{\mathcal{R}}$ as $\mathcal{A}$ - modules we can assume:

1) $Q$ does not depend on the variable $x_3$;

2) $q=x_2^2$;

3) $(1,0,0,0)$ is a double point, which is locally a double cover of the
  plane branched along a curve with a singularity of order at least $4$.

4) $c^2|x_2 z_0^2$;

5) $u=-\frac{w_1z_0}{ec}$.

Moreover $\varphi(E)$ is the line $x_2=x_3=0$.
\end{lemma}

We would like to point out that the coordinates and generators in the
previous lemma can
be chosen such that remark \ref{z-1}, lemma \ref{rc}, remark
\ref{deg2} and lemma
\ref{deg3} will remain valid as it can be traced in the following proof.

{\it Proof.}

First, taking suitable linear combinations of rows and columns of $\alpha$,
we can assume that both $q$ and $Q$ do not depend on the variable
$x_3$, and part 1) is proved.\\

In corollary \ref{rc} we have seen that $vz_{-1}=qz_0$.

As a matter of fact, for every quadric $q'(x_0,x_1,x_2)$ with the property
$$
(**)\ \ \ \ \ \ z_{-1}|q'z_0,
$$
$\frac{q'z_{0}}{z_{-1}}$ is an element
in $\tilde{R}_3$; so if there were two independent quadrics with this
property, we would get two independent elements of 
 ${\tilde{R}}$ in
   degree $3$, and by theorem \ref{FisCM}, for a suitable quadric $q''$
   in the pencil generated by them, we would get a relation
   $g_3(x_i)z_{-1}=q''(x_0,x_1,x_2)z_0$, contradicting corollary
   \ref{resM}. Therefore there can be only one quadric with this property.

We already noticed (cf. remark \ref{cnondividex_0^2z_0}) that $c$
does not divide
$x_0^2z_0$ but divides $x_iz_0$ for every $i \geq 1$; so $q$ is the only
quadric of the form $x_0l(x_1,x_2)+l_1(x_1,x_2)l_2(x_1,x_2)$ such that
$q z_0$ vanishes twice on $E$.

By definition $HE=1$ whence $\varphi (E)$ is a line and there
are two independent linear forms in $\mathbb{P}^3$
vanishing on $E$; in particular there is at least one linear form
in the span of $x_0,x_1,x_2$ vanishing on $E$. Note that $e$ does not divide
$z_0$ (or we could easily find two different quadrics with the
required property, a contradiction).

One of these linear forms belongs to the span of $x_1,x_2$: otherwise
we could assume (up to a change of coordinates) that $x_0$ vanishes on
$E$, and since $q$ is divisible by $e^2$ we get $l_1=l_2=0$ and
$e^2|x_0$, contradicting again the unicity of $q$ (we can take
$x_0x_1$ and $x_0x_2$).

Up to a change of coordinates we can then assume $e|x_2$; $x_2^2$
  fulfills $(**)$, hence $q=x_2^2$, and part $2)$ is proved. \\

Note that since $x_3z_0=-Qz_{-1}$, and $e \! \! \not\vert z_0$ (else
the base point of the canonical system would also be a base point of
the bicanonical system), $e^2|x_3$, and $\varphi(E)=\{x_2=x_3=0\}$.

Let us write $C=C_1+C_2$ (and
accordingly $c=c_1c_2$) where $C_2$ is the greatest common divisor
of $C$ and the divisor of $z_0$, hence obviously $C_1 \neq 0$. Since
$c$ divides $x_i
z_0$ for all $i \geq 1$, and not $x_0z_0$, so
$C_1$ maps to the point $(1,0,0,0)$. \\
By assumption 2) $(1,0,0,0)$ is an isolated singular point of
$\Sigma$, and since the equation of $\Sigma$ is given by the determinant of
$\alpha$, i.e. by $Q^2x_2^2=x_3QG+x_3^2B$, $Q$ is invertible in a
neighbourhood of it. 
Therefore $\Sigma$ has a double point in
$(1,0,0,0)$, which is not a rational double point (otherwise $C =0$). By the form of the
equation of $\Sigma$ we see immediately that the tangent cone has then an
equation of the form
$(x_2+ax_3)^2=0$.

After a linear change of coordinates  we can assume that the tangent
cone is given by $x_2^2=0$. Notice
that this coordinate change "corrupts" the previous choices,
i.e. statement $1)$ and $2)$ do not hold anymore, but we can easily act
on the rows and columns of $\alpha$ in order to "recover" them.

We can now consider the double point as a (local analytic) double cover of the
plane branched on a singular curve with a singularity of order at
least $3$. Assume that the singularity is a triple point. Then by
\cite{BPV} it has to be at
least a $(3,3)$ - point (since otherwise we would have a rational
double point). By our
equations the tangent direction is
$\{x_3=0\}$ and after a blow - up there is again a triple point
exactly on the intersection of the exceptional divisor with the strict
transform of the "tangent" line $\{x_3=0\}$.

In particular the strict trasform of $x_3$ with respect to this blow
up pulls back on $\tilde{S}$ to a divisor with some common component
with $C_1$, while the one of $x_1$ has no common component with
$C_1$. This however contradicts the equality 
$x_3z_0=-Qe^2c$ since $c$ divides  $x_1z_0$ and $Q$ is invertible at
$C$. 
This shows that the branch curve has a singular
point of order at least $4$, and part $3)$ is proven.\\

This in particular implies (looking at the equation of $\Sigma$) that
$c_1^4|x_2^2$, so $c_1^2|x_2$  and $c^2|x_2z_0^2$; this
proves part  $4)$ of the statement.

It is now clear that $\frac{w_1z_0}{ec}$ is a
holomorphic section in $3(H+E)$. Moreover, as we already observed, 
$z_0$ does not vanish on $E$ (or $S$ would have base points for the
bicanonical system), and $w_1$ vanishes on $E$ with multiplicity $1$ 
(it vanishes there because it is multiple of $x_2$, with multiplicity
$1$ or $w_1$ would induce an element of 
$\tilde{R}$). So $\frac{w_1z_0}{ec}$ does not vanish on $E$;
we can therefore choose $u=\frac{w_1z_0}{ec}$, and part $5)$ is
proven.

\hfill $\underline{Q.E.D.}$

\

The choice of $u$ allows us immediately to write the matrix

$$A:=\left(
\begin{array}{cccc}
y_1 & y_2 & y_3 & w_1\\
w_0 & w_1 & Q(y_i) & u
\end{array}
\right),
$$

and to notice that the $2 \times 2$ minors of $A$ are relations in
${\mathcal R}$.

The ring is in fact the canonical ring of a
surface with $7$ generators, so it has codimension $4$. There is no
structure theorem for rings of this codimension, but Reid noticed that
most of them have $9$ relations (joked by $16$ syzygies) that can be
expressed in some ``formats'' (cf. \cite{Rei}, \cite{Rei2})
that help in the study of the deformations.

An important format introduced by Reid is the "rolling factor" format;
we try now to recall shortly how it is defined, referring to the above
quoted papers by Reid for a more detailed treatement and other
examples.

\begin{defi}
One says that a sequence of $9$ equations $f_1, \ldots, f_9$
(usually joked by $16$ syzygies and defining a Gorenstein ring of
codimension $4$, but we do not need this here) is in the
``rolling factor'' format if:

1) $f_1,\ldots f_6$ can be written as the (determinants of the) $2
\times 2$ minors of a $2 \times 4$ matrix $A$;

2) $f_7$ is in the ideal generated by the entries of the first row of
    $A$ (for the matrix $A$ above, it means that $f_7$ can be written as
a linear combination $ay_1+by_2+cy_3+dw_1=0$);

3) $f_8$ is obtained ``rolling'' $f_7$, i.e. taking a linear
combination with the same coefficients, but of the entries of the
second row of $A$ (in our case $f_8$ can be
chosen as $aw_0+bw_1+cQ+du=0$).

4) $f_8$ is in the ideal generated by the entries of the first rows of
    $A$ and $f_9$ is obtained ``rolling'' $f_8$.
\end{defi}

Notice that there can be different ways to ``roll'' the same equation,
but all equivalent up to the equations given by the minors of $A$.

\begin{rem}
If ${\mathcal R}={\mathbb C}[t_0, \ldots t_n]/I$ is an integral
domain, $I$ contains the $2 \times 2$ minors of a $2 \times n$ matrix,
$n \geq 2$, and contains one equation $f$ in the ideal generated by the
entries of the first row of $A$, $I$ must contain also the equation
obtained ``rolling'' $f$.
\end{rem}

This remark will be useful for the next theorem, where we will compute
$f_7$ and prove that $f_7$ can be "rolled" twice, obtaining $f_8$ and $f_9$.

Philosophically, all the three relations come from a single relation in a
bigger ring, that is, the last equation in corollary \ref{rc}.

\begin{theo}\label{canring}
The canonical ring of a surface with $p_g=4$, $K^2=7$, such that
$|K|$ has one simple
base point and $\varphi_{K}$ is birational, is of the form ${\mathcal
R}:={\mathbb
C}[y_0,y_1,y_2,y_3,w_0,w_1,u]/I$ where $I$ is generated by the $2
\times 2$ minors of the
matrix $A$ above, and three more polynomials; one
of degree $4$ of the form $-w_1^2+B(y_i)+\sum \mu_{ijk} y_i y_j w_k$, and the
other $2$ (of respective degrees $5,6$) obtained rolling it twice (so they have the form
$-w_1u+\ldots$, $-u^2+\ldots$).
\end{theo}

In the next section we will write explicitly these equations.

{\it Proof.}

We know already all the generators and the relations in degree $\leq
  3$; moreover we know that all the minors of $A$ are relations of
${\mathcal R}$. An
  easy dimension count shows that there is one relation missing in degree
  $4$; this relation is in fact induced by the rank condition, as follows:

in corollary \ref{rc} we have seen that

$$
vz_0=B z_{-1} -G z_0.
$$

Using the fact that $v=\frac{qz_0}{z_{-1}}=\frac{x_2^2z_0}{e^2c}$
(cf. corollary \ref{rc}
and lemma \ref{geo})

$$
x_2^2 \frac{z^2_0}{e^2c}=B z_{-1} -G z_0.
$$
Multipliying by $\frac{e^2}{c}$
we get the equality
$$
w_1^2 =B(y_i) -\frac{G(y_i)}{ec} z_0.
$$

Using that ${Q^2x_2^2=x_3QG+x_3^2B}$ is singular in $(1,0,0,0)$ (cf.
proof of lemma
\ref{geo}) and recalling that $G=G(x_0,x_1,x_2)$ (cf. remark
\ref{deg2}), we see that the coefficient
of $x_0^3$ in $G$ has to be zero and therefore $\frac{G(y_i)}{ec}
z_0$ can be written as $-\sum \mu_{ijk} y_i y_j w_k$ for
suitable coefficients $\mu_{ijk}$: this provides a non
trivial relation in degree $4$ of the form
$$
(\#\#) \ \ \ \ \ w_1^2 =B(y_i) + \sum \mu_{ijk} y_i y_j w_k
$$
which obviously is  not in the ideal generated by the minors of $A$.

There are no further relations in degree $4$ because
they would force a new generator of ${\mathcal R}$ in degree $4$, which
is excluded by theorem \ref{boundongens}.

We showed in the proof of the last lemma that the singular point
$(1,0,0,0)$ is locally a double cover
of the plane branched along a curve with a singularity of order at least
$4$ and has tangent cone $x_2^2$.

It is easy to verify that this implies that $QG+x_3B$ is contained in the ideal
$(x_2^2,x_2x_1^2,
x_2x_1x_3,x_2x_3^3,x_1^4,x_1x_3^3,x_1^2x_3^2,x_1x_3^3,x_3^4)$, i.e.
the monomials $x_0^3,
x_0^2x_1, x_0^2x_2, x_0x_1^2$ do not appear in $G$, whence $B$ is a quartic in
$\mathbb{P}^3$ such that the monomials
$x_0^4, x_0^3x_1,x_0^3x_3$ have coefficient zero.

Implementing this in $(\#\#)$ we find that the right side can be chosen
(up to adding
some element of the ideal generated by the minors of $A$) to be in
the square of the ideal
generated by the first row of $A$, hence it can be rolled twice.

By theorem \ref{boundongens}, we know that the ideal $I$ of relations
is generated in degree $\leq 6$. We have three elements $f_1,f_2,f_3$
in $I_3$, $f_4,f_5,f_6,f_7$ in $I_4$, $f_8$ in $I_5$, $f_9$ in $I_6$.

Let $I'$ be the ideal $(f_1,\ldots,f_9)$, $R'$ be the quotient ring
${\mathbb C}[y_0,y_1,y_2,y_3,w_0,w_1,u]/I'$. To show $I'=I$ it
suffices to show $I'_k=I_k$ $\forall k \leq 6$, or equivalently, $\dim
R'_k \leq \dim R_k$ for $k \leq 6$. This is a calculation done by
Macaulay 2 (cf. Appendix 1, where the verification is done using  
the equations in theorem \ref{extrasymmetric}).

\hfill $\underline{Q.E.D.}$

\

\section{The canonical ring of surfaces with $K^2 = 7$, $p_g =4$
birational to a sextic: explicit
computations}

Let $S$ be as in the previous section, i.e., a minimal
surface with $K^2=7$, $p_g=4$, such that the canonical system has one
simple base point
and the canonical map is birational.

In the last section we have shown that
the image $\varphi(S)=\Sigma$ of the canonical map has an equation of
the form 

$$
- x_2^2 Q^2 + x_3QG + x_3^2B
$$

where $Q$ is a quadric, $G$ is a cubic, and both do not depend on $x_3$.

Moreover, we have seen in the proof of theorem \ref{canring} that 
the monomials $x_0^3,
x_0^2x_1, x_0^2x_2, x_0x_1^2$ do not appear in $G$, and that $B$ is a quartic in
$\mathbb{P}^3$ such that the monomials
$x_0^4, x_0^3x_1,x_0^3x_3$ have coefficient zero.

In the following three subsections we write the canonical
ring in three explicit ways, using different formats.

\subsection{Rolling factors format}

The proof of theorem \ref{canring} provides immediately the
equations in the "rolling factors" format that we have introduced in the
previous section.

We have defined
$$
A:=
\left(
\begin{array}{cccc}
y_1 & y_2 & y_3 & w_1 \\
w_0 & w_1 & Q(y_0,y_1,y_2) & u
\end{array}
\right).
$$

We write

$$
G = kx_1^3 + x_2Q_0(x_1,x_2) + x_0x_2l(x_1,x_2),
$$
$$
B = x_2C(x_0,x_1,x_2,x_3) + x_3^2Q_3(x_0,x_1,x_3) + x_3x_1Q_2(x_0,x_1) +
x_1^2Q_1(x_0,x_1),
$$
where $k \in \mathbb{C}$, $l$ is linear, the $Q_i$'s are quadratic and $C$
is a cubic. \\

With the above notation, the proof of theorem \ref{canring} gives the
following description for the canonical ring:

\begin{theo}\label{rollfactor}
Let $S$ be a minimal (smooth, connected) surface with $K^2 = 7$ and
$p_g(S) = 4$, such that the
canonical system $|K_S|$ has one (simple) base point
$x \in S$ and the canonical map is birational. Then the canonical
ring $\mathcal{R}(S)$ of $S$ can
be written as

$$
\mathbb{C}[y_0,y_1,y_2,y_3,w_0,w_1,u]/ \mathcal{I},
$$

where $\deg(y_i,w_j,u) = (1,2,3)$ and $\mathcal{I}$ is
generated by the $2 \times 2$ - minors of the matrix

$$
\left(
\begin{array}{cccc}
y_1 & y_2 & y_3 & w_1 \\
w_0 & w_1 & Q(y_0,y_1,y_2) & u
\end{array}
\right),
$$

and by

\begin{multline*}
- w_1^2 + kw_0y_1^2 + w_1Q_0(y_1,y_2) + w_1y_0l(y_1,y_2) +
y_2C(y_0,y_1,y_2,y_3)+\\
+ y_3^2Q_3(y_0,y_1,y_3) + y_1y_3Q_2(y_0,y_1)+ y_1^2Q_1(y_0,y_1),
\end{multline*}

\begin{multline*}
- w_1u + kw_0^2y_1 + uQ_0(y_1,y_2) + uy_0l(y_1,y_2) + w_1C(y_0,y_1,y_2,y_3)+\\
+ y_3Q(y_0,y_1,y_2)Q_3(y_0,y_1,y_3) + w_0y_3Q_2(y_0,y_1) +y_1w_0Q_1(y_0,y_1),
\end{multline*}

\begin{multline*}
- u^2 + kw_0^3 + w_1Q_0(w_0,w_1) + uy_0l(w_0,w_1)
+uC(y_0,y_1,y_2,y_3)+\\
+ Q^2(y_0,y_1,y_2)Q_3(y_0,y_1,y_3)+ w_0Q(y_0,y_1,y_2)Q_2(y_0,y_1)
+w_0^2Q_1(y_0,y_1).
\end{multline*}

\end{theo}

{\it Proof.}
This is just the explicit expression of the computation in theorem
\ref{canring}.

\hfill $\underline{Q.E.D.}$

\

\begin{rem}
1) We remark that the above equations are not exactly "rolled", but
only up to changing
$f_8$ by a suitable combination of the $2 \times 2$ - minors of $A$.
Nevertheless we
prefer to leave the equations in the above form because they are more readable.

2) Our goal is to show that this canonical ring can be deformed to
a canonical ring of surfaces with the same invariants but with base point
free canonical system.

Other ways of writing the same canonical ring might under this aspect
be more convenient.
\end{rem}

\subsection{The $AM(^tA)$ -format.}

Following the notation of M. Reid in \cite{Rei2}, the $AM(^tA)$ -
format is another way to describe
the canonical ring of a Gorenstein variety of codimension $4$
(defined by $9$ relations).It  was introduced by D. Dicks and M. Reid. \\

We briefly recall the definition of the $AM(^tA)$ - format, referring
for details to
\cite{Rei}, \cite{Rei2}, \cite{Rei3}.

\begin{defi}
Assume we are given a commutative ring $R = \mathbb{C}[x_1, \ldots,
x_n]/ (f_1,\ldots,f_9)$. Then
we call this presentation in {\em $AM(^tA)$ - format} iff there is a
$2 \times 4$ - matrix $A$ such
that $f_1, \ldots f_6$ are the $2 \times 2$ - minors of $A$ and there is
a symmetric $4 \times  4$ -
matrix $M$ such that
$$
AM(^tA) = \left(
\begin{array}{cc}
f_7 & f_8 \\
f_8 & f_9
\end{array}
\right).
$$
\end{defi}

This format is slightly more restricted than the previous one. In fact, the
equations given in the $AM(^tA)$ format are automatically in the
"rolling factor"
format, whereas equations given in the "rolling factor" format
can be expressed in an $AM(^tA)$ format if and only if the $3$
supplementary relations are
"quadratic forms in the rows of $A$". This is equivalent to saying
that the first one
(the one we denoted by $f_7$) is in the square of the ideal generated
by the entries
of the first row of $A$, what happens to occur in our specific situation.\\

\bigskip
We are now ready to state our result.

\begin{theo}\label{AMtA}
Let $S$ be a minimal (smooth, connected) surface with $K^2 = 7$ and
$p_g(S) = 4$, such that
the canonical system $|K_S|$ has one (simple) base point
$x \in S$ and the canonical map is birational. Then the canonical
ring $\mathcal{R}$ of $S$ is of the form ${\mathbb
C}[y_0,y_1,y_2,y_3,w_0,w_1,u]/{\mathcal I}$ (deg
$(y_i,w_j,u)=(1,2,3)$) and can be presented
in the $AM(^tA)$ - format with

$$
A = \left(
\begin{array}{cccc}
y_1 & y_2 & y_3 & w_1 \\
w_0 & w_1 & Q(y_0,y_1,y_2) & u
\end{array}
\right),
$$

and

$$
M = \left(
\begin{array}{cccc}
m_{11} & m_{12} & m_{13} & k_1y_0 \\
  & m_{22} & m_{23} & k_2y_0 \\
  &  & m_{33} & L \\
sym &  &  & -1
\end{array}
\right).
$$

Here $k_1, k_2 \in \mathbb{C}$, $L = L(y_0,y_1,y_2,y_3)$ is linear
and $m_{ij}$ are quadratic forms
depending on the variables $y_0,y_1,y_2,y_3, w_0,w_1$ as follows:

$$
\begin{array}{lll}
m_{11} & = & m_{11}(y_0,y_1,w_0,w_1), \\
m_{12} & = & a_1w_1 + a_2 y_0y_1 + a_3y_0y_2 + a_4y_1^2, \\
m_{13} & = &m_{13}(y_0,y_1), \\
m_{22} & = & b_1w_1 + b_2y_0y_2 + b_3y_1^2 + b_4y_1y_2 + b_5 y_2^2, \\
m_{23} & = & m_{23}(y_0,y_1,y_2,y_3), \\
m_{33} & = & m_{33}(y_0,y_1,y_3),
\end{array}
$$

where $a_i, b_j \in \mathbb{C}$.
\end{theo}

As we will see in the proof, the coefficients of the entries of $M$
can be explicitly determined from
the equation of $\Sigma$. \\

{\em Proof.}
First we note that $A$ was already found in theorem \ref{rollfactor}.
In order to
write down the matrix $M$,
we have to set up some more notation. \\

We write explicitly (keeping the terminology which was introduced earlier)

$$
Q_0(y_1,y_2) = q_{11}y_1^2 + q_{12}y_1y_2 + q_{22}y_2^2,
$$

$$
l(y_1,y_2) = l_1y_1 + l_2y_2,
$$
where $q_{ij}, l_k \in \mathbb{C}$. \\

We remark that the coefficient of $y_0^2$ in $Q$ has to be different from
zero because of the assumption that
the generalized tacnode $o$ is not contained in $\Gamma$. Therefore
we can write

\begin{multline*}
C(y_0,y_1,y_2,y_3) = l'(y_0,y_1,y_2,y_3)Q + y_0(Q'_{11}y_1^2 +
Q'_{12}y_1y_2 + Q'_{22}y_2^2) + \\
y_0y_3l''(y_1,y_2,y_3) + sy_1^3 + y_2Q''(y_1,y_2) + y_3Q'''(y_1,y_2,y_3),
\end{multline*}

where $s, Q'_{ij} \in \mathbb{C}$, $l', l''$ are linear forms and
$Q'', Q'''$ are quadratic forms. \\

Now, a rather lengthy, but straightforward calculation shows that for $M=$

$$
\left(
\begin{array}{cccc}
kw_0 + Q_1 + q_{11}w_1 & \frac{1}{2} (q_{12}w_1 + y_0(Q'_{11}y_1 +
Q'_{12}y_2) + sy_1^2) &
\frac{1}{2} Q_2 & \frac{1}{2} l_1y_0 \\
   & q_{22}w_1 + Q'_{22}y_0y_2 + Q'' & \frac{1}{2} (y_0l'' + Q''') &
\frac{1}{2} l_2y_0\\
   &  & Q_3 & \frac{1}{2} l' \\
sym & & & -1
\end{array}
\right),
$$

the conditions rank $A \leq 1$ and $AM(^tA) = 0$ define
the ideal $\mathcal{I}$ of theorem $\ref{rollfactor}$. 

\hfill $\underline{Q.E.D.}$

\

\begin{rem}
In \cite{Rei2} M. Reid shows that, given a polynomial ring ${\mathbb
C}[x_0, \ldots, x_n]$
and a Gorenstein ${\mathbb C}[x_0, \ldots,
x_n]-$algebra ${\mathcal R}$ of codimension $4$ and 
presented in
$AM(^tA)$ - format, all the syzygies are induced by $A$ and $M$. In
particular in this
case the $A M ( ^t A)
$-format is {\it flexible}, i.e. every deformation of the matrices $A$ and $M$
preserving the symmetry of $M$ induces a flat deformation of ${\mathcal R}$.\\
In order to find a flat family ${\mathcal S}/T$ of surfaces such that
$S_{t_0}$ is a
surface with $K^2=7$, $p_g=4$ such that the canonical system has one
base point and the
canonical map  is birational, whereas  $S_t$ is a surface such that
the canonical system has no base points for every $t \neq t_0$, it
would be sufficient
to find a deformation of the above matrices, which induces the right ring.\\
Unfortunately this does not work, since by \cite{CatBu} the canonical
ring of $S_t$ (for
$t \neq t_0$) is generated in degrees $1$ and $2$. In particular, if
$\mathcal{R}_t$ is such a deformation of $\mathcal{R}$, for $t \neq
0$, one of the three relations
in degree three has to eliminate the generator $u$. But, considering
only deformations of
$\mathcal{R}$ induced by deformations of $A$ and $M$, the relations
in degree three of
$\mathcal{R}_t$ are given by the two by two minors of a matrix
${A}_t$ with the following
degrees in the entries:

$$
\left(
\begin{array}{ccc}
1 & 1 & 1\\
2 & 2 & 2
\end{array}
\right),
$$

which obviously cannot eliminate an element $u$ having degree $3$.
\end{rem}

\subsection{The antisymmetric and extrasymmetric format.}

The aim of this section is to give a third description of our
canonical ring $\mathcal{R} =
\mathbb{C}[y_0,y_1,y_2,y_3,w_0,w_1,u]/ \mathcal{I}$.

\begin{defi}
1) A $6 \times 6$ - matrix $P$ is called {\em antisymmetric and
extrasymmetric} if and only if it has the form 
$$
P =
\left(
\begin{array}{cccccc}
  0 &  a & b   & c & d & e \\
    &  0 & f   & g & h & d \\
    &    & 0   & i & pg & pc \\
    &    &     & 0 & qf & qb \\
    &    &     &   & 0 & pqa \\
-sym&    &     &   &   & 0
\end{array}
\right).
$$

2) Assume we are given a commutative ring ${\mathcal R} =
\mathbb{C}[x_1, \ldots,
x_n]/{\mathcal I}$. If there is an antisymmetric and extrasymmetric matrix
$P$ whose $4
\times 4$ pfaffians
generate the ideal ${\mathcal I}$, we say that the ring has an {\em
antisymmetric and extrasymmetric format
given by $P$}.
\end{defi}

This format was first introduced by D. Dicks and M. Reid in a less general form
\cite{Rei}. The above more general form appeared in \cite{Rei3}.

An easy computation shows that the ideal of the $15$ ($4 \times 4$)
pfaffians of an antysimmetric and extrasymmetric
matrix is in fact generated by nine of them.

Under suitable generality assumptions (e.g. that ${\mathcal R}$ be
Gorenstein of codimension $4$ and that the $9$ equations be independent),
it should be easy (but rather lengthy) to show,
following the lines of Reid's  argument used in \cite{Rei}
to prove a special case, 
that the format is
flexible (i.e., that a deformation of the entries of $P$ induces a flat
deformation of the ring).

Since however we are only interested in our particular case, it does
not pay off here to perform this
calculation. Hence we will proceed as follows: \\ 
we shall put our canonical ring in an antisymmetric and 
extrasymmetric format, deform the associated
$6
\times 6$ antisymmetric and extrasymmetric matrix (preserving the
extrasymmetry) and
then verify the flatness
of the induced family 
by the constancy of the Hilbert polynomial.

Using the same notation as in the previous section we obtain

\begin{theo}\label{extrasymmetric}
Let $S$ be a minimal (smooth, connected) surface with $K^2 = 7$ and
$p_g(S) = 4$, whose
canonical system $|K_S|$ has one (simple) base point
$x \in S$ and yields a birational canonical map.
Then the canonical ring of $S$ can be presented as $\mathcal{R} =
\mathbb{C}[y_0,y_1,y_2,y_3,w_0,w_1,u]/
\mathcal{I}$, where deg $(y_i,w_j,u) = (1,2,3)$
and  the ideal of relations
$\mathcal{I}$ of $\mathcal{R}$ is generated by the
$4 \times 4$ - pfaffians of the following
antysimmetric and extrasymmetric matrix

$$
P =
\left(
\begin{array}{cccccc}
0 & 0 & w_0 & Q(y_0,y_1,y_2) & w_1 & u \\
  & 0 & y_1 & y_3 & y_2 & w_1 \\
  &  & 0 & -u+\overline{C} & y_3\overline{Q}_3(w_0,y_0,y_1,y_3) & Q\overline{Q}_3 \\
  & &  & 0 & y_1\overline{Q}_1(w_0,w_1,y_0,y_1,y_3) & w_0\overline{Q}_1 \\
& &  &  & 0 & 0 \\
-sym &  &  &  & & 0
\end{array}
\right),
$$

where $Q, \overline{Q}_1, \overline{Q}_3$ are quadratic
forms of a subset of the given variables as indicated,
and $\overline{C}$ is a cubic form. Moreover $\overline{C}$ does not
depend on $u$ and the  $w_iy_j$'s for $j
\leq 1$, $\overline{Q}_1$ does not depend on $y_3^2$.
\end{theo}

{\em Proof.}

Recall that the coefficient of $x_0^2$ in $Q$ is different from
zero. This allows us to choose $s \in {\mathbb C}$ and linear forms
$\overline{l}_1(y_0,y_1)$, $\overline{l}_2(y_0,y_1,y_2)$ such that

$$
Q_2(y_0,y_1) = sQ(y_0,y_1,y_2) + y_1\overline{l}_1(y_0,y_1) +
y_2\overline{l}_2(y_0,y_1,y_2).
$$
Therefore we can write

$$
\overline{C} = q_{12}w_0y_2 + q_{22}w_1y_2 + y_0l(w_0,w_1) + C(y_0,y_1,y_2,y_3) +
y_1y_3\overline{l}_2(y_0,y_1,y_2);
$$

$$
\overline{Q}_1 = kw_0 + q_{11}w_1 + Q_1(y_0,y_1)+y_3\overline{l}_1(y_0,y_1)
$$

$$
  \overline{Q}_3 = - Q_3 - sw_0.
$$

Now the rest of the proof is a straightforward calculation.

\hfill $\underline{Q.E.D.}$

\

\section{An explicit family.}

In this section we will find an explicit deformation of the canonical
ring $\mathcal{R}$ to the
canonical ring $\mathcal{S}$ of a surface with $K^2 = 7$, $p_g = 4$,
such that the canonical system
has no base points. \\
Before doing this we have to recall some of the results (\cite{CatBu}) on
surfaces with $K^2 = 7$, $p_g = 4$, whose canonical system is base
point free.

\bigskip

Let $X$ be a nonsingular surface with $K^2 = 7$, $p_g = 4$, such that
the canonical system is base
point free, hence the canonical map is a birational morphism onto a
septic surface in
${\mathbb P}^3$. We denote by
$\mathcal{S}$ the canonical ring of
$X$. As in the previous case, we denote by $y_i$ an appropriate basis of
$H^0(X,{\mathcal O}_X(K))$.
We set ${\mathcal A}:= {\mathbb C}[y_0,y_1,y_2,y_3]$.

\begin{theo}\label{fabrizio}
1)
$\mathcal{S}$ has a minimal resolution
as $\mathcal{A}$ - module given by the matrix

$$
\alpha = \left(
\begin{array}{ccc}
d_1d_2y_0 + (d_3d_4 + d_2^2)y_1 + (d_2d_3 + d_1d_4)y_2 & d_4y_1 &
d_1y_0 + d_2y_1 + d_3y_2 \\
d_4y_1 & y_0 & y_2 \\
d_1y_0 + d_2y_1 + d_3y_2 & y_2 & y_1
\end{array}
\right),
$$

where $d_1,d_2,d_3,d_4$ are arbitrary quadratic forms in $y_i$. \\
2) $\alpha$ satisfies the rank condition $\Lambda^2(\alpha) =
\Lambda^2(\alpha')$, where $\alpha'$ is obtained by deleting the
first row of $\alpha$, and
therefore induces a unique ring structure on $\mathcal{S}$ as
quotient of $\mathcal{B}=\mathcal{A}[w_0,w_1]$ by
the three relations given by
$$
\alpha
\left(
\begin{array}{c}
1 \\
w_0 \\
w_1
\end{array}
\right) = 0
$$

and three more relations expressing $w_iw_j$ as linear combination of
the other monomials whose
coefficients are determined by the adjoint matrix of $\alpha$. \\
3) The surfaces with $K^2 = 7$, $p_g = 4$ such that the canonical
system is base point free form an irreducible unirational component of dimension
$36$ in the moduli space $\mathfrak{M}_{K^2 = 7, p_g = 4}$ of surfaces
with $K^2 = 7$, $p_g = 4$. 
\end{theo}

For a more precise formulation we refer to the original articles
\cite{CatBu} or \cite{Homalg}.

\begin{rem}
1) The previous theorem implies in particular that the canonical ring
$\mathcal{S}$  is Gorenstein in codimension
$3$. Hence by the classical result of Buchsbaum and Eisenbud (cf.
\cite{BuEi}) the ideal
$\mathcal{I} \subset \mathcal{A}[w_0,w_1] =: \mathcal{B}$ defining
$\mathcal{S}$ can be minimally
generated by the
$2k \times 2k$ - Pfaffians of a skewsymmetric $(2k + 1) \times (2k +
1)$ - matrix. Writing down
explicitly the above six defining equations, we can see that
$\mathcal{I}$ is generated by only five
of them, hence in our case $k=2$. \\
2) More precisely, we have a selfdual resolution of $\mathcal{S}$ as
$\mathcal{B}$ - module as
follows

$$
0 \rightarrow \mathcal{B}(-9)  \stackrel{f_2}{\rightarrow}
\mathcal{B}(-5)^3 \oplus
\mathcal{B}(-6)^2
\stackrel{f_1}{\rightarrow} \mathcal{B}(-4)^3 \oplus
\mathcal{B}(-3)^2 \stackrel{f_0}{\rightarrow}
\mathcal{B}
\rightarrow
\mathcal{S} \rightarrow 0,
$$

where $f_1$ is alternating, $f_0$ is given by the Pfaffians of $f_1$
and $f_2 =\ ^tf_0$. \\
3) Vice versa, assume we have a $\mathbb{C}$ - algebra $\mathcal{S}$,
which admits a resolution as above: then under suitable open condition
$\mathcal{S}$ is the
canonical ring of a surface $X$ with $p_g = 4$, $K^2 = 7$ and free
canonical system. 
\end{rem}

Our aim is now to take the matrix $P$ in antisymmetric and extrasymmetric
format and try to
find a deformation $P_t$ of $P$
with the following properties: \\

0) For $t\neq 0$ in at least one of the Pfaffians of degree $3$ the
generator $u$ appears with a non
zero coefficient; \\

1) for $t\neq 0$ there is a skewsymmetric $5 \times 5$ -matrix $Q_t$,
such that the
ideal ${\mathcal J}_t$ generated by the
$4 \times 4$ - Pfaffians of $Q_t$ coincides with $\mathcal{I}_t \cap
\mathcal{B}$, where
$\mathcal{I}_t$ is the ideal generated by the
$4 \times 4$ - Pfaffians of $P_t$, for every $t \neq 0$; \\

2) the entries of $Q_t$ have the right degrees, i.e. $Q_t$ defines a map

$$
\mathcal{B}(-5)^3 \oplus
\mathcal{B}(-6)^2
\stackrel{Q_t}{\rightarrow} \mathcal{B}(-4)^3 \oplus \mathcal{B}(-3)^2;
$$

3) $\mathcal{S}_t := {\mathbb C}[y_0,y_1,y_2,y_3,w_0,w_1,u]/
\mathcal{I}_t$ is a flat
family.

\bigskip
 From the previous remark we conclude that $\mathcal{S}_t$ (for $t
\neq 0$), constructed as above,
is the canonical ring of a surface $X$ with
$p_g = 4$, $K^2 = 7$ and free canonical system. Hence once we have
found a deformation $P_t$ as
above, we have explicitly deformed the surfaces with 
$p_g = 4$, $K^2 = 7$, such that $|K_S|$ has one base point and
induces a birational map to the
surfaces with $p_g = 4$, $K^2 = 7$ and free canonical system.

\bigskip
The natural attempt now is to deform the entries of $P$ preserving
the extrasymmetry (this should give automatically the flatness), so
that one of the relations in degree $3$ eliminates $u$. At first
glance we see a natural way: putting the deformation parameter $t$ in the
only entry of degree $0$ (replacing the $0$ therein) and replacing the symmetrical ``zero'' to preserve the
extrasymmetry. This
actually works and we have the following

\begin{theo}\label{famiglia}
Let $P$ be an antisymmetric and extrasymmetric matrix as in theorem
\ref{extrasymmetric}.
Consider the 1-parameter
family of rings
${\mathcal S}_t={\mathbb C}[y_0,y_1,y_2,y_3,w_0.w_1,u]/{\mathcal
I}_t$ where the ideal
${\mathcal I}_t$ is given by the $4 \times 4$ pfaffians  of the
antysimmetric and extrasymmetric matrix
$$
P_t =P+
\left(
\begin{array}{cccccc}
0 & t &0 & 0 & 0 & 0 \\
  & 0 & 0 & 0 & 0 & 0 \\
  &  & 0 & 0 & 0 & 0 \\
  &  &  & 0 & 0 & 0 \\
  &  & &  & 0 & t\overline{Q}_1\overline{Q}_3 \\
  -sym &  &  &  & & 0
\end{array}
\right).
$$

This is a flat family and  describes a flat deformation of the
surface corresponding to the matrix $P$ to surfaces with $p_g=4$, $K^2=7$ and
with $|K|$ base point free.
For $t \neq 0$,
${\mathcal S}_t$ is isomorphic to ${\mathbb
C}[y_0,y_1,y_2,y_3,w_0.w_1]/J_t$
where $J_t$ is the ideal generated by the $4 \times 4$
pfaffians of  the
matrix
$$
\left(
\begin{array}{ccccc}
0 & y_1 & -y_3 & -t^2w_0 & -tQ \\
   & 0 & y_2 & t^3\overline{Q}_3 & tw_1 \\
   &  & 0 & t^2w_1 & t^2\overline{Q}_1  \\
   &  &  & 0 & -t^3c-t^2y_1Q+t^2y_3w_0\\
-sym &  &  &  & 0
\end{array}
\right).
$$
\end{theo}

{\it Proof.}

We have to check that the properties 0)-3) described above are
fullfilled by our pair
$P_t, Q_t$.

0): We immediately see that the $4 \times 4$ - Pfaffian obtained
eliminating the last two rows and columns eliminates the generator $u$;

1): once we eliminate $u$ it is easy to write down explicitely
$\mathcal{I}_t \cap \mathcal{B}$ for $t \neq 0$ and check that it
coincides with the given ideal $J_t$ (cf. Appendix 2);

2): obvious;

3): we expect that flatness holds in general once we preserve the
extrasymmetry. 

In our particular case, flatness follows auitomatically since:

a) $X_t:=$Proj$({\mathcal B}/J_t)$ has dimension $\leq 2$ by
semicontinuity. 

b) Then $Q_t$ and its pfaffians give a resolution of $J_t$ by
\cite{BuEi}.

c) Therefore the Hilbert polynomial of $J_t$ is the same as the one of
$J_0$, then the family is flat.

Finally, since the property that $X_t$ has only R.D.P.'s as
singularities is open, $X_t$ is the canonical model of a surface of
general type as required.

\hfill $\underline{Q.E.D.}$

\

\begin{rem}
As we have already said in the introduction, this particular
degeneration was studied by Enriques in his book \cite{Enr}.
There Enriques states that such a degeneration should exist and suggest
a way to construct such a family degenerating the canonical
images. With the help of Macaulay 2 we have explicitely computed the
degenerations of the canonical images corresponding to our family and
we have found that the degeneration is not the one ``predicted'' by
Enriques. In the following we recall briefly 
(see \cite{Enr} for the details) Enriques' prediction and point out
where is the difference.
\end{rem}

The canonical image of a surface with $K^2=7$, $p_g=4$ with base point
free canonical system and  birational canonical morphism 
is a surface in ${\mathbb P}^3$ of degree $7$ with a singular curve of
degree $7$ and genus $4$ having a triple point. Moreover the
adjoint quadric is a quadric cone (with vertex in the singular point
of the curve) whose intersection with the surface is given by the
above curve of degree $7$ counted twice. 
 
As we have already seen, the canonical image of a surface with
$K^2=7$, $p_g=4$ with one simple case point for the canonical system
and canonical map birational
is a surface in ${\mathbb P}^3$ of degree $6$ with a singular curve of
degree $2$ and a generalized tacnode. In this case there is an adjoint 
plane that is the plane through the conic ($\{x_3=0\}$ in our notation).

Enriques suggests to add to the sextic the adjoint plane. 
The intersection of this plane
and the sextic is given by the singular conic and the line image of
$E$ both counted twice;
this gives a reducible septic with a triple conic and a double line of 
``tacnodal type'' (i.e. near a general point of the line the surface
has two branches tangent on the line).

Enriques states that it is possible to construct a family of septics
with a singular curve of degree $7$ and genus $4$ having a triple
point that degenerates to the above configuration so that the singular
septic degenerates to the union of the conic (counted three times) and
the line.

We wrote the canonical images of the family described by the $4 \times
4$ pfaffians of the $5 \times 5$ skewsymmetric matrix in theorem 
\ref{famiglia}, and, as anticipated, we did not get the situation
predicted by Enriques. In fact the family of septics degenerates to the
union of the sextic canonical image of the limit surface with a plane,
but instead of the plane predicted by Enriques we have gotten the
plane $\{x_2=0\}$ (the reduced tangent cone of the tacnode).

In fact, it is quite easy to compute also the degeneration of the
adjoint quadric: consider the resolution of the canonical ring of a
surface with base point free canonical system as ${\mathcal A}-$module
as in theorem \ref{fabrizio}. From this resolution one can immediately
see that the adjoint quadric must be the determinant of the
right-bottom $2$ by $2$ minor of the resolution matrix in theorem
\ref{fabrizio} 
(i.e. $y_0y_1-y_2^2$ in the coordinates chosen there). This minor
depends only on the two relations in degree $3$ of the canonical ring.

It is now easy to compute it for our family: we have just to write
down the two relations in degree $3$, write the $2$ by $2$ matrix of
the coefficients of $w_0$ and $w_1$ in this two equations, and then
compute the determinant.

In the notation of theorem \ref{famiglia} the two relations in degree $3$
(for $t \neq 0$) can be written as:
$$
y_1w_1-w_0y_2+ty_3\overline{Q}_3;
$$
$$
ty_1\overline{Q}_1+w_1y_3-Qy_2.
$$
The equation of the quadric cone depends clearly on the coefficients
of $w_0$ and $w_1$ in $\overline{Q}_1$ and $\overline{Q}_1$ (but notice
that it is in every case independent of $y_0$, so it cannot be a smooth
quadric, as expected), but the two equations degenerate respectively to
$y_1w_1-w_0y_2
$ and
$
w_1y_3-Qy_2
$; 
the $2$ by $2$ minor degenerates then to
$$
\begin{pmatrix}
-y_2 & y_1 \\
0 & y_3
\end{pmatrix}
$$
and the quadric cone degenerate to the union of the adjoint plane and
the tacnodal plane.

Geometrically we could say that the septic degenerates to the union of
a sextic and a plane, the adjoint quadric to the union of the same plane and
a different plane (the tacnodal plane); the two ``identical'' planes
``simplify'' and we are left with the sextic and his adjoint plane.

\

\

Address of the authors until August 31, 2001:

Mathematisches Institut der Universit\"at G\"ottingen \\
Bunsenstr. 3-5 \\
37073 G\"ottingen\\

Address of the authors after September 1, 2001:

Mathematisches Institut der Universit\"at Bayreuth \\
Lehrstuhl Mathematik VIII\\
Universit\"atstr. 30 \\
95447 Bayreuth\\

E-mail addresses:\\
Ingrid.Bauer@uni-bayreuth.de\\
Fabrizio.Catanese@uni-bayreuth.de\\
Roberto.Pignatelli@uni-bayreuth.de

\newpage

{\small
\begin{verbatim}


-- APPENDIX 1
-- This script checks that the relations we found in Theorem 2.13 are all the
-- relations of the canonical ring till degree 6, so they are all the
-- relations by theorem 1.8.
-- First we write the ring: with all the variables and parameters we need
R=QQ[u,w_1,w_0,y_0,y_1,y_2,y_3,
     a0,a1,a2,b0,d0,d1,t,Q,Q1,Q3,c,q00,q01,q02,q11,q12,q22,
     q100,q101,q103,q111,q113,q3,q300,q301,q303,q311,q313,q333,
     c000,c001,c002,c003,c011,c012,c013,c022,c023,c033,
     c111,c112,c113,c122,c123,c133,c222,c223,c333,
     MonomialOrder=>Lex
      ]
-- Now we write the matrix in theorem 3.7
M=matrix{
{0,0,w_0,Q,w_1,u},
{0,0,y_1,y_3,y_2,w_1},
{-w_0,-y_1,0,-u+(c+a0*w_1*y_0-a1*w_1*y_1+a2*w_1*y_2+b0*w_0*y_0),y_3*Q3,Q*Q3},
{-Q,-y_3,u-(c+a0*w_1*y_0-a1*w_1*y_1+a2*w_1*y_2+b0*w_0*y_0),0,
               (Q1+d0*w_0+d1*w_1)*y_1,w_0*(Q1+d0*w_0+d1*w_1)},
{-w_1,-y_2,-y_3*Q3,-(Q1+d0*w_0+d1*w_1)*y_1,0,0},
{-u,-w_1,-Q*Q3,-w_0*(Q1+d0*w_0+d1*w_1),0,0}
};
pfaff=pfaffians(4,M);
-- Here we restrict to 9 pfaffians and check that they are enough to
-- generate the whole pfaffian ideal (the second line gives ``true''
-- as output)
pfaff9=submatrix(gens(pfaff),,{0,1,5,2,6,9,3,4,7});
gens(pfaff) % ideal(pfaff9)==0
-- Then we write explicitly all the polynomials in the matrix
ourideal:=substitute(pfaff9,{
          Q=>y_0*y_0+q01*y_0*y_1+q02*y_0*y_2+q11*y_1*y_1+q12*y_1*y_2+q22*y_2*y_2,
          Q1=>q100*y_0*y_0+q101*y_0*y_1+q103*y_0*y_3+q111*y_1*y_1+q113*y_1*y_3,
          Q3=>q3*w_0+q300*y_0*y_0+q301*y_0*y_1+q303*y_0*y_3+q311*y_1*y_1+
              q313*y_1*y_3+q333*y_3*y_3,
          c=>c000*y_0*y_0*y_0+c001*y_0*y_0*y_1+c002*y_0*y_0*y_2+c003*y_0*y_0*y_3+
             c011*y_0*y_1*y_1+c012*y_0*y_1*y_2+c013*y_0*y_1*y_3+c022*y_0*y_2*y_2+
             c023*y_0*y_2*y_3+c033*y_0*y_3*y_3+c111*y_1*y_1*y_1+c112*y_1*y_1*y_2+
             c113*y_1*y_1*y_3+c122*y_1*y_2*y_2+c123*y_1*y_2*y_3+c133*y_1*y_3*y_3
             +c222*y_2*y_2*y_2+c223*y_2*y_2*y_3+c333*y_3*y_3*y_3
          })
-- Here we define the ideal of the monomials in all the degrees till 6 
linear:=ideal(y_0,y_1,y_2,y_3);
quadrics:=ideal(w_0,w_1)+linear^2;
cubics:=ideal(mingens(ideal(u)+linear*quadrics));
quartics:=ideal(mingens(quadrics^2+linear*cubics));
quintics:=ideal(mingens(quadrics*cubics+linear*quartics));
sextics:=ideal(mingens(cubics*cubics+quadrics*quartics+linear*quintics));
-- finally we compute a system of generators, degree by degree, of the
-- resulting quotient. All of them turn out to be composed exactly by
-- P_n elements (resp. 4,12,26,47,75), that concludes the argument in the 
-- proof of theorem 2.13 
K=mingens ideal((gens linear) % gb ourideal);
twoK=mingens ideal((gens quadrics) % gb ourideal);
threeK=mingens ideal((gens cubics) % gb ourideal);
fourK=mingens ideal((gens quartics) % gb ourideal);
fiveK=mingens ideal((gens quintics) % gb ourideal);
sixK=mingens ideal((gens sextics) % gb ourideal);

restart

-- APPENDIX 2

R=QQ[t,u,y_1..y_3,Q,Q1,Q3,w_0..w_1,c,
      Degrees=>{1,3,1,1,1,2,2,2,2,2,3}]

-- the 4 x 4 - Pfaffians of the matrix M define the relations
-- of the canonical ring of a surface with K^2 = 7, p_g = 4 such that the
-- canonical map has one base point and is birational.


M=matrix{
{0,0,w_0,Q,w_1,u},
{0,0,y_1,y_3,y_2,w_1},
{-w_0,-y_1,0,-u+c,y_3*Q3,Q*Q3},
{-Q,-y_3,u-c,0,Q1*y_1,w_0*Q1},
{-w_1,-y_2,-y_3*Q3,-Q1*y_1,0,0},
{-u,-w_1,-Q*Q3,-w_0*Q1,0,0}};
M=map(R^{-1,-2,3:0,1},R^{-2,-1,3:-3,-4},M)

-- we calculate the 4 x 4 - Pfaffians of M and extract the 9
-- ``important'' as above;
pfaff15=pfaffians(4,M);
pfaff9=submatrix(gens(pfaff15),,{0,1,5,2,6,9,3,4,7});
gens(pfaff15) % ideal(pfaff9)==0

-- we write the sixteen syzygies of them
syzs=syz pfaff9

-- we define the matrix $Mt = M + tM1$, which is the deformation of $M$ whose
-- pfaffians we want to understand;

M1=matrix(R,{
          {0,1,0,0,0,0},
          {-1,0,0,0,0,0},
          {0,0,0,0,0,0},
          {0,0,0,0,0,0},
          {0,0,0,0,0,Q3*Q1},
          {0,0,0,0,-Q3*Q1,0}
          });
M1=map(R^{-1,-2,3:0,1},R^{-2,-1,3:-3,-4},M1)
Mt=M+t*M1;
-- we calculate the 15 Pfaffians of Mt and verify that the same nine 
-- Pfaffians of Pt again generate the whole ideal (the output of the third 
-- line below is ``true'';
defpfaff=pfaffians(4,Mt);
defpfaff9=submatrix(gens(defpfaff),,{0,1,5,2,6,9,3,4,7});
gens(defpfaff) % ideal(defpfaff9)==0

-- If t in different from 0, one can eliminate the variable u using the first 
-- equation in defpfaff9.

elimu=defpfaff9_(0,0)

-- we will use the following trick to eliminate u:
-- u appears only in degrees smaller than 2 in defpfaff9; we
-- multiply defpfaff9 by t^2, and reduce by elimu;
-- what we get is the same ideal defpfaff as before for every t
-- different from zero!
-- Finally we divide by t, wherever it is possible;
-- here it is crucial that we choose a monomial order such that tu
-- is the leading term of elimu: this forces the result to be independent of u;

defpfaff9withoutu:= divideByVariable((t^2*defpfaff9) % elimu,t);

-- the following 5 generators are enough; we choose a strange order in order
-- to get a nicer result;
fiveequationsfortnotzero=submatrix(defpfaff9withoutu,,{5,4,2,3,1});
defpfaff9withoutu % ideal(fiveequationsfortnotzero)==0

-- now we look for the 5x5 matrix inducing these equations as pfaffians:

lookforQ=syz fiveequationsfortnotzero;


-- among the 36 syzygies (Macaulay found a lot of them because he is 
-- considering also the case t=0) one can easily find something that looks
-- interesting
almostQ=submatrix(lookforQ,,{1,3,2,5,9})
-- this matrix is not (yet) antisymmetric; we change coordinates in the source
-- and in the target in order to make it antisymmetric;
one:=matrix(R,{{1}})
diag1=one++one++one++(t*one)++(-1*one)
diag2=one++one++(-1*one)++(-t^2*one)++(-t*one)
Qt=diag1*almostQ*diag2

-- finally we check, whether for t different from 0, the Pfaffians of Qt and
-- the Pfaffians of Mt after having eliminated u generate the same ideal;
pfaffQt=pfaffians(4,Qt);
pfaffQred=divideByVariable(gens pfaffQt,t);
fiveequationsfortnotzero %  pfaffQred==0
pfaffQred % fiveequationsfortnotzero==0
\end{verbatim}
}

\end{document}